\documentclass[a4paper,11pt,reqno]{amsart}
 
\usepackage[latin1]{inputenc}
\usepackage{amsmath,amssymb,amsthm}
\usepackage[margin=2.5cm]{geometry}
\usepackage{enumerate}
\usepackage{hyperref}
\usepackage{xcolor}
\usepackage{tikz}
\usepackage{mathdots}
\usepackage{FloPPPackage}

\newtheorem{thm}{Theorem}[section]

\newtheorem{cor}[thm]{Corollary}
\newtheorem{prop}[thm]{Proposition}

\newtheorem{lem}[thm]{Lemma}

\DeclareMathOperator{\sgn}{sgn}
\DeclareMathOperator{\even}{even}
\DeclareMathOperator{\ASM}{ASM}
\DeclareMathOperator{\Id}{Id}
\DeclareMathOperator{\CT}{CT}
\newcommand{\Sgrp}{\mathfrak{S}}
\newcommand{\AS}{\mathcal{AS}}
\newcommand{\Z}{\mathbb{Z}}
\newcommand{\N}{\mathbb{N}}
\newcommand{\Q}{\mathbb{Q}}
\newcommand{\x}{\mathbf{x}}
\newcommand{\overDelta}[1]{\overline{\Delta}_{\textit{ }#1}}
\newcommand{\underDelta}[1]{\underline{\Delta}_{\textit{ }#1}}

\newcommand{\Tiling}{
	\PlanePartitionWhite{{4,4,3,1},{4,3,0,-1},{4,2},{3,2},{3,1},{2}}
	\draw [fill= white] (0,0) -- ({-2*cos(30)},-1) -- (0,-2);
	\TopBoundary{2}{-4}{-2}{3}
	\TopBoundary{2}{-5}{2}{3}
	\TopBoundary{2}{-6}{1}{3}	
	\TopBoundary{3}{-4}{-2}{3}
	\TopBoundary{3}{-5}{1}{3}
	\TopBoundary{3}{-6}{1}{3}
	\TopBoundary{4}{-2}{-1}{3}
	\TopBoundary{4}{-3}{-2}{3}
	\TopBoundary{4}{-4}{-2}{3}
	\TopBoundary{4}{-5}{1}{3}
	\TopBoundary{5}{-5}{1}{3}	
	\LeftBoundary{2}{-5}{2}{3}
	\LeftBoundary{2}{-6}{1}{3}
	\LeftBoundary{3}{-4}{2}{3}
	\LeftBoundary{3}{-6}{1}{3}
	\LeftBoundary{4}{-4}{2}{3}
	\LeftBoundary{4}{-5}{1}{3}
	\LeftBoundary{5}{-4}{2}{3}
	\LeftBoundary{5}{-5}{1}{3}
	\LeftBoundary{6}{-4}{1}{3}
	\LeftBoundary{3}{0}{4}{3}
	\LeftBoundary{4}{0}{4}{3}
	\LeftBoundary{4}{0}{3}{3}
	\LeftBoundary{4}{0}{2}{3}
	\RightBoundary{0}{-4}{4}{3}
	\RightBoundary{0}{-5}{4}{3}
	\RightBoundary{0}{-6}{4}{3}
	\RightBoundary{0}{-6}{3}{3}
	\RightBoundary{2}{-5}{2}{3}
	\RightBoundary{6}{-3}{3}{3}
	\RightBoundary{6}{-3}{2}{3}
	\RightBoundary{6}{-3}{1}{3}
	\RightBoundary{6}{-4}{1}{3}
}

\title[A new determinant for the $Q$-enumeration of ASMs]{A new determinant for the $Q$-enumeration of alternating sign matrices}
\author{Florian Aigner}
\address{Florian Aigner, Universit\"at Wien, Fakult\"at f\"ur Mathematik, Oskar-Morgenstern-Platz~1, 1090 Wien, Austria}
\email{florian.aigner@univie.ac.at}
\thanks{Supported by the Austrian Science Fund FWF, START grant Y463, SFB grant F50 and J 4387.}

\begin{document}

%%%%%%%%%%%%%%%%%%%%%%%%%%%%%%%%%%%%%%%%%%%%%%%%%%%%%%%%%%%%%%%%%%%%%%%%%%%%%%%%%%%%%%%%%%%%%
%%%%%%%%%%%%%%%%%%%%%%%%%%%%%%%%%%%%%%%%%  Abstract  %%%%%%%%%%%%%%%%%%%%%%%%%%%%%%%%%%%%%%%%
%%%%%%%%%%%%%%%%%%%%%%%%%%%%%%%%%%%%%%%%%%%%%%%%%%%%%%%%%%%%%%%%%%%%%%%%%%%%%%%%%%%%%%%%%%%%%

\begin{abstract}
Fischer provided a new type of binomial determinant for the number of alternating sign matrices involving the third root of unity. In this paper we prove that her formula, when replacing the third root of unity by an indeterminate $q$, actually gives the $(q^{-1}+2+q)$-enumeration of alternating sign matrices.
By evaluating a generalisation of this determinant we are able to reprove a conjecture of Mills, Robbins and Rumsey stating that the $Q$-enumeration is a product of two polynomials in $Q$. Further we provide a closed product formula for the  generalised determinant in the $0$-, $1$-, $2$- and $3$-enumeration case, leading to new proofs of the $1$-, $2$- and $3$-enumeration of alternating sign matrices, and a factorisation in the $4$-enumeration case.
Finally we relate the $1$-enumeration case of our generalised determinant to the determinant evaluations of Ciucu, Eisenk\"olbl, Krattenthaler and Zare, which count weighted cyclically symmetric lozenge tilings of a hexagon with a triangular hole and are a generalisation of a famous result by Andrews. As a result, we obtain alternative proofs of their determinantal evaluations using the Desnanot-Jacobi identity (Dodgson condensation).
\end{abstract}

\maketitle

%%%%%%%%%%%%%%%%%%%%%%%%%%%%%%%%%%%%%%%%%%%%%%%%%%%%%%%%%%%%%%%%%%%%%%%%%%%%%%%%%%%%%%%%%%%%%
%%%%%%%%%%%%%%%%%%%%%%%%%%%%%%%%%%%%%%%%%%  Intro  %%%%%%%%%%%%%%%%%%%%%%%%%%%%%%%%%%%%%%%%%% 
%%%%%%%%%%%%%%%%%%%%%%%%%%%%%%%%%%%%%%%%%%%%%%%%%%%%%%%%%%%%%%%%%%%%%%%%%%%%%%%%%%%%%%%%%%%%%

\section{Introduction}

An \emph{alternating sign matrix} (\emph{ASM}) of size $n$ is an $n \times n$ matrix with entries $-1,0$ or $1$ such that all row and column sums are equal to $1$ and the non-zero entries alternate in each row and column, see Figure \ref{fig: ASM to monotone triangle} (left). These matrices were introduced by Robbins and Rumsey \cite{RobbinsRumsey86} and arose from generalising the determinant to the so-called $\lambda$-determinant. Together with Mills \cite{MillsRobbinsRumsey82} they conjectured that the number of ASMs of size $n$ is
\[
\prod_{i=0}^{n-1} \frac{(3i+1)!}{(n+i)!}.
\]
This was remarkable since there already existed a different family of combinatorial objects, \emph{descending plane partitions} (or DPPs), which were proved by Andrews \cite{Andrews79} to have the same enumeration formula, but which are of a very different nature. In \cite{MillsRobbinsRumsey86} Mills, Robbins and Rumsey conjectured this very enumeration formula for a third family of combinatorial objects, namely \emph{totally symmetric self complementary plane partitions} (or \emph{TSSCPP}s), which was proved by Andrews \cite{Andrews94}.\\

The first proof of the ASM Theorem was found by Zeilberger \cite{Zeilberger96}, who related a constant term formula for ASMs to a constant term formula for TSSCPPs. Shortly thereafter, Kuperberg \cite{Kuperberg96} presented a proof by using a different approach. He used the fact that there exists an easy bijection between ASMs and configurations of the six-vertex model, a model in statistical mechanics. This technique was later used as a standard method to prove enumeration formulas of various symmetry classes and refinements of ASMs.
In \cite{MillsRobbinsRumsey83} Mills, Robbins and Rumsey also proved a formula for the $2$-enumeration and conjectured one for the 3-enumeration of ASMs, where the $Q$-enumeration of ASMs is a weighted enumeration of ASMs with each ASM having the weight of $Q$ to the power of the number of its $-1$'s. Kuperberg \cite{Kuperberg96} was able to prove with his methods the $3$-enumeration and also that the $Q$-enumeration of ASMs is a product of two polynomials in $Q$ which was also conjectured in \cite[Conjecture 4]{MillsRobbinsRumsey83}.\\

Roughly a decade later, Fischer \cite{Fischer07} provided a new proof of the ASM Theorem which was independent of the first two methods. The proof is based on the bijection between ASMs and monotone triangles as well as the \emph{operator formula} \cite[Theorem 1]{Fischer06} which enumerates monotone triangles.
In 2016 Fischer \cite{Fischer16} presented a concise version of her original proof which, except for relying on the \emph{Lindstr\"om-Gessel-Viennot} Theorem at one point, is self-contained.\\

A \emph{monotone triangle} with $n$ rows is a triangular array $(a_{i,j})_{1\leq j \leq i \leq n}$ of integers of the following form,
\[
\begin{array}{c c c c c c c c c c c}
&&&&& a_{1,1}\\
&&&& a_{2,1} && a_{2,2}\\
&&& \iddots && \cdots &&\ddots\\
&& \iddots && \iddots && \ddots &&\ddots\\
& a_{n-1,1} && a_{n-1,2} && \cdots  && \cdots  && a_{n-1,n-1}\\
a_{n,1} && a_{n,2} && a_{n,3} && \cdots  && \cdots && a_{n,n}
\end{array}
\]
such that the entries are weakly increasing along northeast and southeast diagonals, i.e., $a_{i+1,j}  \leq a_{i,j} \leq a_{i+1,j+1}$, and strictly increasing along rows. For an example see Figure \ref{fig: ASM to monotone triangle} (right).
Given an ASM, the $i$-th row from the top of its corresponding monotone triangle records the columns of the ASM with a positive partial column sum of the top $i$ rows. It is easy to see that this yields a bijection between ASMs of size $n$ and monotone triangles with bottom row $1,2,\ldots, n$. We assign to a monotone triangle a weight $Q^\Sigma$ where $\Sigma$ is the number of entries $a_{i,j}$ of the triangle such that $a_{i+1,j} < a_{i,j} < a_{i+1,j+1}$. For ASMs this weight is exactly $Q$ to the power of the number of $-1$'s.
\begin{figure}
\[
\begin{pmatrix}
0 & 0 & 0 & 1 & 0\\
0 & 1 & 0 & 0 & 0\\
1 & -1 & 1 & 0 & 0\\
0 & 1 & 0 & -1 & 1\\
0 & 0 & 0 & 1 & 0\\
\end{pmatrix}
 \qquad \leftrightarrow \qquad
 \begin{array}{c cc cc cc cc}
&&&& 4\\
&&& 2 && 4\\
&& 1 && 3 && 4\\
& 1 && 2 && 3 && 5\\
1 && 2 && 3 && 4 && 5
 \end{array}
\]
\caption{An ASM of size $5$ and its corresponding monotone triangle.}
\label{fig: ASM to monotone triangle}
\end{figure}
\newline

The core of this paper is the following theorem.

\begin{thm}
\label{thm: main thm}
Let $k$ be an integer and $n$ be a positive integer and define
\begin{equation*}
\label{eq: det in betrachtung}
d_{n,k}(x,q) :=\det_{1 \leq i,j \leq n}\left(\binom{x+i+j-2}{j-1}\frac{1-(-q)^{j-i+k}}{1+q}\right).
\end{equation*}
Then the $(q^{-1}+2+q)$-enumeration of alternating sign matrices is equal to $d_{n,1}(0,q)$.
\end{thm}

The above determinant appeared first in \cite[p. 559]{Fischer18} for $k=1$ and was shown by Fischer to count the number of ASMs for $x=0, k=1$ and $q$ being a primitive third root of unity. By introducing the variable $x$, which was suggested in \cite{Fischer18}, and the integer parameter $k$ in the determinant, we are able to write $d_{n,k}(x,q)$ as a closed product formula for arbitrary $k$ and $q$ being a primitive second root of unity (Theorem \ref{thm: second root}), primitive third root of unity (Theorem \ref{thm: third root}), primitive fourth root of unity (Theorem \ref{thm: fourth root}) or primitive sixth root of unity (Theorem \ref{thm: sixth root}), which was conjectured for $k=1$ in \cite{Fischer18}. For $q=1$ we provide in Theorem \ref{thm: first root} a factorisation of the determinant as a polynomial in $x$.
 Compared to other known determinantal formulas for the $Q$-enumeration of ASMs for which the evaluation is rather complicated, the evaluation of the determinant considered in this paper turns out to be very easy and thus leads immediately to the known formulas for the $0$-, $1$-, $2$- and $3$-enumeration of alternating sign matrices by setting $x=0$.\\

In Theorem \ref{thm: structural thm} we prove a general factorisation result which states that the determinant $d_{n,k}(x,q)$ factors for arbitrary $q$ into a power of $q$, a polynomial $p_{n,k}(x) \in \mathbb{Q}[x]$ which factorises into linear factors and a polynomial $f_{n,k}(x,q) \in \mathbb{Q}[x,q]$ which is given recursively.
For $k=1$, Theorem \ref{thm: structural thm} implies that the determinant $d_{n,k}(x,q)$ can be written as a product of two Laurent polynomials in $q$ with coefficients in $\mathbb{Q}[x]$, which was conjectured in \cite{Fischer18}.
As a direct consequence we obtain that the generating function of ASMs with respect to the number of $-1$'s is a product of two polynomials in $Q=q^{-1}+2+q$.\\

Surprisingly the determinant $d_{n,k}(x,q)$ is connected to the famous determinant by Andrews \cite{Andrews79}  and its generalisation by Ciucu, Eisenk\"olbl, Krattenthaler and Zare \cite{CiucuEisenkoelblKrattenthalerZare01} in the following way:
\begin{equation}
\label{eq: det verbindung}
d_{n,3-k}(x,q^2)=q^{-n}  \det_{1 \leq i,j \leq n} \left(\binom{x+i+j-2}{j-1}+q^k \delta_{i,j}\right),
\end{equation}
where $q$ is a primitive sixth root of unity. This fact was first conjectured for $k=2,4$ in \cite{Fischer18} and is remarkable because of two reasons. Firstly, for $q$ being a primitive sixth root of unity and $k$ an integer the evaluations of the determinants have a closed product formula. This fact is very easy to prove for the left hand side of \eqref{eq: det verbindung} by using the Desnanot-Jacobi identity which is also known as the \emph{Condensation method}. For the right hand side this method is however not applicable and the proof relies on the \emph{method of identification of factors}. Secondly, the determinant by Ciucu et al. is a weighted enumeration of cyclically symmetric  lozenge tilings of hexagons with a triangular hole, where $q$ corresponds to the weight and $x$ is the side length of the hole. In \cite{Krattenthaler06} a bijection between these tilings and $d$-DPPs, ``one parameter generalisations'' of DPPs introduced by Andrews \cite{Andrews79, Andrews80}, was established.
We prove \eqref{eq: det verbindung} by connecting $d_{n,k}(x,q)$ to a determinantal formula for the enumeration of column strict shifted plane partitions (CSSPPs), which are in bijection with $d$-DPPs and using the determinantal formula $\det(M_{\textnormal{ASM}})$ for the weighted enumeration of ASMs by Behrend, Di Francesco and Zinn-Justin \cite[Proposition 1]{BehrendDiFrancescoZinnJustin12}. As a consequence we obtain an alternative proof of Theorem \ref{thm: main thm} which is based on the six-vertex model approach as well as an alternative proof for the factorisation of the right hand side of \eqref{eq: det verbindung} by the Desnanot-Jacobi identity.

Very recently Fischer proved in \cite{Fischer19b} that this determinant also enumerates alternating sign trapezoids, which are one parameter generalisations of ASTs. Hence the determinant $d_{n,k}(x,q)$ suggests a one parameter refinement for ASMs and might be of help in finding one.
\\

The paper is structured in the following way.
In Section \ref{sec: herleitung} we follow the steps of \cite{Fischer18, FischerRiegler15} and prove Theorem \ref{thm: main thm}. Section \ref{sec: q allgemein} contains a description of the factorisation of the determinant $d_{n,k}(x,q)$ for general $q$ leading to a proof that the $(q^{-1}+2+q)$-enumeration of ASMs is a product of two polynomials. In Section \ref{sec: andere wurzeln} we present  product formulas for $d_{n,k}(x,q)$ for $q$ being a primitive second, third, fourth or sixth root of unity and present a factorisation for $q=1$. This leads to new proofs for the $1$-, $2$- and $3$-enumeration of ASMs. 
Finally in Section \ref{sec: other dets} we relate the determinant $d_{n,k}(x,q)$ to the weighted enumeration of CSSPPs. This allows us to connect $d_{n,k}(x,q)$ with another determinantal formula for ASMs described in \cite{BehrendDiFrancescoZinnJustin12}, the Andrews determinant and its generalisation by Ciucu, Eisenk\"olbl, Krattenthaler and Zare. The paper ends with an appendix containing a list of specialisations of $d_{n,k}(x,q)$ which turn out to be known enumeration formulas.

%Finally in Section \ref{sec: verbindungen} we relate the determinant $d_{n,k}(x,q)$ to the Andrews determinant and its generalisation by Ciucu, Eisenk\"olbl, Krattenthaler and Zare. The paper ends with an appendix containing a list of specialisations of $d_{n,k}(x,q)$ which turn out to be known enumeration formulas.

\section{A determinantal formula for the number of ASMs using the operator formula}
\label{sec: herleitung}

%%%%%%%%%%%%%%%%%%%%%%%%%%%%%%%%%%%%%%%%%%%%%%%%%%%%%%%%%%%%%%%%%%%%%%%%%%%%%%%%%%%%%%%%%%%%%
%%%%%%%%%%%%%%%%%%%%%%%%%%%%%%%%%%%%%%%%  CT of ASMs %%%%%%%%%%%%%%%%%%%%%%%%%%%%%%%%%%%%%%%% 
%%%%%%%%%%%%%%%%%%%%%%%%%%%%%%%%%%%%%%%%%%%%%%%%%%%%%%%%%%%%%%%%%%%%%%%%%%%%%%%%%%%%%%%%%%%%%

Let $f$ be a function in one or $n$ variables respectively. For the rest of this paper we fix the following notations,
\begin{align*}
E_x (f)(x):=& f(x+1) \qquad &\textit{shift operator}, \\
\overDelta{x}:=&E_x -\Id \qquad &\textit{forward difference},\\
\underDelta{x}:=& \Id - E_x^{-1} \qquad &\textit{backward difference}.
\end{align*}
We denote by $\AS_{x_1,\ldots,x_n}f(x_1,\ldots,x_n)$ the antisymmetriser of $f$ with respect to $x_1,\ldots,x_n$ which is defined as
\[
\AS_{x_1,\ldots,x_n}f(x_1,\ldots,x_n) = \sum_{\sigma \in \Sgrp_n} \sgn(\sigma) f(x_{\sigma(1)}, \ldots, x_{\sigma(n)}).
\]
The constant term of a function $f$ with respect to the variables $x_1, \ldots, x_n$ is denoted by $\CT_{x_1, \ldots, x_n} \left( f \right)$ and $[x^i]f(x)$ denotes the coefficient of $x^i$ in $f(x)$.
Further we are using at various points the multi-index notation, i.\,e., a bold variable $\mathbf{x}$ refers to a vector $ \mathbf{x}= (x_1, \ldots, x_n)$ and $\mathbf{x}^{\mathbf{a}}$ is defined as $\mathbf{x}^{\mathbf{a}}:=\prod_{i=1}^nx_i^{a_i}$.\\

In this section we deduce the determinantal expression for the $Q$-enumeration of ASMs. It is based on results of \cite{FischerRiegler15} and generalises results in \cite{Fischer18}.
The starting point is a weighted version of the operator formula for monotone triangles.

\begin{thm}[{\cite[Theorem 1]{Fischer10}}]
\label{thm: monotone triangles}
 The generating function of monotone triangles with bottom row $\mathbf{k}=(k_1,\ldots, k_n)$ with respect to the $Q$-weight is given by evaluating the polynomial
 \begin{equation*}
 \label{eq: operator formula}
  M_n(x_1,\ldots,x_n):= \prod_{1 \leq i < j \leq n}(Q \Id+ (Q-1)\overDelta{x_i}+\overDelta{x_j}+\overDelta{x_i}\overDelta{x_j})
  \prod_{1 \leq i < j \leq n} \frac{x_j-x_i}{j-i},
 \end{equation*}
at $\mathbf{x}= \mathbf{k}$.
\end{thm}

Following \cite{FischerRiegler15}, we can rewrite the operator formula as a constant term expression. We want to point out, that the original statement \cite[Proposition 10.1]{FischerRiegler15} is obtained by replacing $x_i$ by $1+z_i$ and $Q$ by $X$ in the next proposition. 
\begin{prop}[{\cite[Proposition 10.1]{FischerRiegler15}}]
\label{prop: constant term formula}
The number of monotone triangles with bottom row $\mathbf{k}=(k_1,\ldots, k_n)$ with respect to the $Q$-weight is the constant term of
\[
\AS_{x_1,\ldots,x_n} \left( \prod_{i=1}^n (1+x_i)^{k_i} \prod_{1 \leq i < j \leq n}(Q+(Q-1)x_i+x_j +x_i x_j) \right)
\prod_{1 \leq i < j \leq n}(x_j-x_i)^{-1}.
\]
\end{prop}

%%%%%%%%%%%%%%%%%%%%%%%%%%%%%%%%%%%%%%%%%%%%%%%%%%%%%%%%%%%%%%%%%%%%%%%%%%%%%%%%%%%%%%%%%%%%%
%%%%%%%%%%%%%%%%%%%%%%%%%%%%%%%%%%%%%%%%%  CT to Det %%%%%%%%%%%%%%%%%%%%%%%%%%%%%%%%%%%%%%%% 
%%%%%%%%%%%%%%%%%%%%%%%%%%%%%%%%%%%%%%%%%%%%%%%%%%%%%%%%%%%%%%%%%%%%%%%%%%%%%%%%%%%%%%%%%%%%%

The following theorem allows us to transfer the antisymmetriser in the above proposition into a determinant; it is a variation by Fischer of Equation (D.3) in \cite{FonsecaZinnJustin08}.

\begin{thm}[{\cite[Theorem 17]{Fischer18}}]
\label{thm: antisym}
Let $f(x,y)= qx-q^{-1}y$ and $h(x,y)=x-y$. Then 
\begin{equation*}
\AS_{w_1, \ldots,w_n} \frac{\prod\limits_{1 \leq i < j \leq n}f(w_i,w_j)}{\prod\limits_{1 \leq i \leq j \leq n} h(w_j,y_i)f(w_i,y_j)}
= \frac{\det\limits_{1 \leq i,j \leq n}\left(\frac{1}{f(w_i,y_j)h(w_i,y_j)} \right)}{\prod\limits_{1 \leq i<j \leq n}h(y_j,y_i)}.
\end{equation*}
\end{thm}

By setting $w_i = \frac{x_i+1+q^{-1}}{x_i+1+q}$, $Q=q^{-1}+2+q$ and taking the limit of $y_i \rightarrow 1$ for all $1 \leq i \leq n$, Theorem \ref{thm: antisym} becomes
\begin{multline*}
\frac{(-1)^{\frac{n(n+1)}{2}}}{(q-q^{-1})^{\frac{n(n+3)}{2}}}
 \AS_{x_1,\ldots,x_n}
\left( \frac{\prod\limits_{1 \leq i < j \leq n} (Q+(Q-1)x_i+x_j +x_ix_j)\prod\limits_{i=1}^n (x_i+1+q)^2}{\prod\limits_{i=1}^n (1+x_i)^{n+1-i}} \right)\\
= q^n \lim_{y_1,\ldots,y_n \rightarrow 1} \det_{1 \leq i,j \leq n}\left(
\frac{1}{\left(y_j - \frac{x_i+1+q^{-1}}{x_i+1+q} \right)\left(y_j - q^2\frac{x_i+1+q^{-1}}{x_i+1+q} \right)} \right)\prod_{1 \leq i<j \leq n}(y_j-y_i)^{-1}.
\end{multline*}
Hence, by Proposition \ref{prop: constant term formula} with $k_i=i$ for each $i$, the $(q^{-1}+2+q)$-enumeration of ASMs is given by
\begin{multline}
\label{eq: asm as det 1}
\CT_{x_1, \ldots, x_n}\left((-1)^{\frac{n(n+1)}{2}}q^n(q-q^{-1})^{\frac{n(n+3)}{2}} \prod_{i=1}^n (1+x_i)^{n+1}(x_i+1+q)^{-2}\right.\\
\left.\times \lim_{y_1,\ldots,y_n \rightarrow 1} \det_{1 \leq i,j \leq n}\left(
\frac{1}{\left(y_j - \frac{x_i+1+q^{-1}}{x_i+1+q} \right)\left(y_j - q^2\frac{x_i+1+q^{-1}}{x_i+1+q} \right)} \right)
\prod_{1 \leq i < j \leq n} (x_j-x_i)^{-1}(y_j-y_i)^{-1}\right).
\end{multline}
Using the partial fraction decomposition $\frac{1}{(y-a)(y-b)}=\frac{1}{a-b}\left( \frac{1}{y-a}-\frac{1}{y-b}\right)$ we can rewrite the determinant in \eqref{eq: asm as det 1} as
\begin{multline*}
\label{eq: det anders ausgedrueckt}
(1-q^2)^{-n}\prod_{i=1}^n (x_i+1+q)^2(x_i+1+q^{-1})^{-1}\\
\det_{1 \leq i,j \leq n}\left(
\frac{1}{y_j(x_i+1+q) - (x_i+1+q^{-1})}-\frac{1}{y_j(x_i+1+q) - q^2(x_i+1+q^{-1})} \right).
\end{multline*}
The following lemma allows us to evaluate the limit in \eqref{eq: asm as det 1}.\\

\begin{lem}[{\cite[Eq. (43)--(47)]{BehrendDiFrancescoZinnJustin12}}]
\label{lem: det at limit}
Let $f(x,y)=\sum\limits_{i,j \geq 0}c_{i,j}x^iy^j$ be a formal power series in $x$ and $y$. Then
\[
\lim_{\substack{x_1,\ldots,x_n \rightarrow x\\y_1,\ldots,y_n \rightarrow y}}
\frac{\det\limits_{1 \leq i,j \leq n}(f(x_i,y_j))}{\prod\limits_{1 \leq i < j \leq n} (x_j-x_i)(y_j-y_i)}
= \det_{0 \leq i,j \leq n-1}\left( [u^iv^j]f(x+u,y+v) \right).
\]
%where $[u^iv^j]f(x+u,y+v)$ denotes the coefficient of $u^iv^j$ in $f(x+u,y+v)$.
\end{lem}

By the above lemma (and since $\CT_{x_1, \ldots, x_n}$ is simply $\lim_{x_1,\ldots,x_n \rightarrow 0}$) we need to calculate the coefficient of $x^iy^j$ in
\begin{equation}
\label{eq: ausdruck in det}
\frac{1}{(y+1)(x+1+q) - (x+1+q^{-1})}-\frac{1}{(y+1)(x+1+q) - q^2(x+1+q^{-1})}.
\end{equation}
Using the geometric series expansion in $x$ and $y$ of \eqref{eq: ausdruck in det}, we obtain for the coefficient of $x^iy^j$ 
\begin{equation*}
\label{eq: cij}
\binom{j}{i}(1+q)^{-i-1}(q-1)^{-j-1}(-1)^jq^{j+1}-\sum_{k=i}^{i+j}(-1)^k\binom{k}{j}\binom{j}{i-k+j}\frac{(1+q)^{k-i-j-1}}{(1-q)^{j+1}}.
\end{equation*}
Putting the above together, it follows that the $(q^{-1}+2+q)$-enumeration of ASMs of size $n$ is given by
\begin{equation*}
\label{eq: asm as det 2}
(1+q)^{-n}(-q)^{-\binom{n}{2}} \det_{0 \leq i,j \leq n-1}
\left( \binom{j}{i} (-1)^jq^{j+1}+ \sum_{k=0}^{n-1} \binom{k+i}{j}\binom{j}{k}(-1-q)^{k+i-j} \right).
\end{equation*}
Finally we use the following identity which is due to Fischer \cite[p. 599]{Fischer18} and can be proved using basic properties of the binomial coefficient
\begin{multline*}
\binom{j}{i} (-1)^jq^{j+1}+ \sum_{k=0}^{n-1} \binom{k+i}{j}\binom{j}{k}(-1-q)^{k+i-j}\\
= \sum_{k=0}^{n-1}(-1)^i\binom{i}{k}\binom{-k-1}{j}(q^{j+1}(-1)^k+q^k(-1)^j).
\end{multline*}
Hence the $(q^{-1}+2+q)$-enumeration of ASMs is
\begin{multline}
\label{eq: asm as det 3}
\det \left(
\left( \binom{i}{j}(-1)^{i+j} \right)_{0 \leq,i,j \leq n-1} \times
\left( \binom{i+j}{j}\frac{1-(-q)^{j-i+1}}{1+q} \right)_{0 \leq i,j \leq n-1} \right)\\
= \det_{1 \leq i,j \leq n} \left( \binom{i+j-2}{j-1}\frac{1-(-q)^{j-i+1}}{1+q} \right),
\end{multline}
(where simple properties such as $\binom{-k-1}{j}=(-1)^j\binom{k+j}{j}$, and the fact that $\left( \binom{i}{j}(-1)^{i+j} \right)_{0 \leq i,j \leq n-1}$ is lower triangular with only $1$'s on the diagonal, have been used), which proves Theorem \ref{thm: main thm}.

%%%%%%%%%%%%%%%%%%%%%%%%%%%%%%%%%%%%%%%%%%%%%%%%%%%%%%%%%%%%%%%%%%%%%%%%%%%%%%%%%%%%%%%%%%%%%
%%%%%%%%%%%%%%%%%%%%%%%%%%%%%%%%%%%%  general structure  %%%%%%%%%%%%%%%%%%%%%%%%%%%%%%%%%%%% 
%%%%%%%%%%%%%%%%%%%%%%%%%%%%%%%%%%%%%%%%%%%%%%%%%%%%%%%%%%%%%%%%%%%%%%%%%%%%%%%%%%%%%%%%%%%%%

\section{A generalised $(q^{-1}+2+q)$-enumeration}
\label{sec: q allgemein}

We introduce to the determinant in \eqref{eq: asm as det 3} a variable $x$, as suggested in \cite{Fischer18}, and further a parameter $k \in \Z$. The determinant of our interest is then $d_{n,k}(x,q):=\det\left( D_{n,k}(x,q) \right)$, where $D_{n,k}(x,q)$ is defined by 
\[
D_{n,k}(x,q):=\left(\binom{x+i+j-2}{j-1}\frac{1-(-q)^{j-i+k}}{1+q}\right)_{1 \leq i,j \leq n}.
\]
The evaluation of this determinant will generally follow two steps. The first step is to guess a formula for $d_{n,k}(x,q)$ using a computer algebra system. The second step is to use induction and the Desnanot-Jacobi Theorem, which is sometimes also called the \emph{condensation method}.

\begin{thm}[Desnanot-Jacobi]
\label{thm: Condensation method}
Let $n$ be a positive integer, $A$ an $n \times n$ matrix and denote by $A_{j_1,\cdots, j_k}^{i_1,\cdots, i_k}$ the submatrix of $A$ in which the $i_1,\cdots, i_k$-th rows and $j_1,\cdots, j_k$-th columns are omitted.
Then holds
 \begin{equation*}
 \label{eq: condensation}
 \det A \det A_{1,n}^{1,n} = \det A_1^1 \det A_n^n -\det A_1^n \det A_n^1.
 \end{equation*}
\end{thm}

Deleting the first and/or last row and the first and/or last column of the matrix $D_{n,k}(x,q)$, and taking determinants, gives the following expressions.
\begin{align}
\label{eq: row column deleting}
 \det(D_{n,k}(x,q)_1^1)&=\det(D_{n-1,k}(x+2,q))\binom{x+n}{n-1},\\
 \nonumber \det(D_{n,k}(x,q)_n^n)&=\det(D_{n-1,k}(x,q)),\\
 \nonumber \det(D_{n,k}(x,q)_{1,n}^{1,n})&=\det(D_{n-2,k}(x+2,q))\binom{x+n-1}{n-2},\\
 \nonumber \det(D_{n,k}(x,q)_n^1)&=\det(D_{n-1,k-1}(x+1,q)),\\
 \nonumber \det(D_{n,k}(x,q)_1^n)&=\det(D_{n-1,k+1}(x+1,q))\binom{x+n-1}{n-1}.
 \end{align}
We will prove it in the case of $D_{n,k}(x,q)_1^1$. The other cases are proved analogously.
By definition we have
\begin{multline*}
\det(D_{n,k}(x,q)_1^1) = \det_{2 \leq i,j \leq n}\left(\binom{x+i+j-2}{j-1}\frac{1-(-q)^{j-i+k}}{1+q}\right)\\
= \det_{1 \leq i,j \leq n-1}\left(\binom{x+i+j}{j}\frac{1-(-q)^{j-i+k}}{1+q}\right).
\end{multline*}
 Since the $i$-th row is divisible by $(x+i+1)$ we can factor it out for all $1 \leq i \leq n-1$. Further we factor out $j^{-1}$ from the $j$-th column for $1 \leq j \leq n-1$ and obtain
\begin{multline*}
\det(D_{n,k}(x,q)_1^1) = \prod_{i=1}^{n-1} \frac{x+i+1}{i}
\det_{1 \leq i,j \leq n-1}\left(\binom{x+i+j}{j-1}\frac{1-(-q)^{j-i+k}}{1+q}\right)\\
=\binom{x+n}{n-1}\det(D_{n-1,k}(x+2,q)).
\end{multline*}
By applying \eqref{eq: row column deleting} to Theorem \ref{thm: Condensation method} we obtain
\begin{multline}
 \label{eq: condensed equation}
 (n-1)d_{n,k}(x,q)d_{n-2,k}(x+2,q)\\
=(x+n)d_{n-1,k}(x,q)d_{n-1,k}(x+2,q)-(x+1)d_{n-1,k+1}(x+1,q)d_{n-1,k-1}(x+1,q).
 \end{multline}
We use at certain points of this section the $-q$ analog $[j]_{-q}$ of an integer $j \in \mathbb{Z}$, which is defined by 
\begin{equation}
\label{eq: -q-analog}
[j]_{-q} =\frac{1-(-q)^j}{1+q}=
\begin{cases}
\sum\limits_{i=0}^{j-1} (-q)^i \qquad & j >0,\\
0 & j=0,\\
-\sum\limits_{i=j}^{-1}(-q)^i & j<0.
\end{cases}
\end{equation}
The entries of the matrix $D_{n,k}(x,q)$ are polynomials in $x$ and Laurent polynomials in $q$ and hence the same is true for the determinant $d_{n,k}(x,q)$, i.e., $d_{n,k}(x,q) \in \Q[q,q^{-1},x]$.

\begin{lem}
\label{lem: divisibility by p(n,k)}
Let $n,k$ be positive integers. Then the determinant $d_{n,k}(x,q)$ is as an element in $\Q[q,q^{-1},x]$ divisible by 
\begin{align*}
 \prod_{l=0}^{\left\lfloor\frac{n-k-1}{2}\right\rfloor} (x+k+2l+1).
\end{align*}
\end{lem}

\begin{proof}
Let $l$ be a non-negative integer. The $(i,j)$-th entry of $D_{n,k}(x,q)$ is divisible by $(x+k+2l+1)$
for $i \leq k+2l+1 \leq i+j-2$, because of the binomial coefficient, and for $i=j+k$ since $[j+k-i]_{-q}=[0]_{-q}=0$. Now choose $l$ with $0 \leq l \leq \left\lfloor\frac{n-k-1}{2}\right\rfloor$ and set  $x=-k-2l-1$. Then the $(i,j)$-th entry of $D_{n,k}(-k-2l-1,q)$ is equal to $0$ for all $k+l+1 \leq i \leq k+2l+1\leq n$ and $j \geq l+1$. Therefore the $(k+l+1)$-st up to the $(k+2l+1)$-st row are linearly dependent and the determinant is henceforth equal to $0$.
\end{proof}

\begin{lem}
\label{lem: transponierte matrix}
Let $n,k$ be non-negative integers. Then the following identity holds
\[
d_{n,-k}(x,q)= (-1)^{n(k-1)}q^{-n k}d_{n,k}(x,q).
\]
\end{lem}

\begin{proof}
Let $\sigma \in S_n$ be a permutation. Then 
\begin{align*}
\prod_{i=1}^n\binom{x+i+\sigma(i)-2}{\sigma(i)-1} =
 \prod_{i=1}^n\binom{x+\sigma(i)+i-2}{i-1}.
\end{align*}
This implies that we can replace the binomial coefficient $\binom{x+i+j-2}{j-1}$ in the determinant $d_{n,k}(x,q)$ with $\binom{x+i+j-2}{i-1}$ without changing the determinant.
By using this fact and factoring out the factor $-(-q)^j$ from the $j$-th column and $(-q)^{-i-k}$ from the $i$-th row for all rows and columns, we obtain
\begin{multline*}
d_{n,-k}(x,q)= \det_{1 \leq i,j \leq n} \left(\binom{x+i+j-2}{j-1}\frac{1-(-q)^{j-i-k}}{1+q}\right)\\
= \prod_{i=1}^n(-q)^{-i-k} \prod_{j=1}^n(-1)(-q)^{j} 
\det_{1 \leq i,j \leq n} \left(\binom{x+i+j-2}{i-1}\frac{-(-q)^{i+k-j}+1}{1+q}\right) \\
=(-1)^{n(k-1)}q^{-n k} \det \left(D_{n,k}(x,q)^T\right)=(-1)^{n(k-1)}q^{-n k} d_{n,k}(x,q).
\end{multline*}
\end{proof}

\begin{cor}
\label{cor: k=0 and odd n}
Let $n$ be an odd positive integer. Then $d_{n,0}(x,q) = 0$.
\end{cor}

With the above two lemmas at hand we can prove the following structural theorem.

\begin{thm}
\label{thm: structural thm}
The determinant $d_{n,k}(x,q)$ has the form 
\begin{equation}
\label{eq: general form}
d_{n,k}(x,q)= q^{c(n,k)}p_{n,k}(x)f_{n,k}(x,q),
\end{equation}
with
\begin{align*}
p_{n,k}(x)&= \prod_{i=1}^{n-1} \frac{\left\lfloor \frac{i}{2}\right\rfloor!}{i!} \prod_{i=0}^{\left\lfloor \frac{n-|k|-1}{2}\right\rfloor}(x+|k|+2i+1),\\
c(n,k)&= \begin{cases}
0  \qquad & k>0, n\leq k,\\
n\, k  \qquad & k<0, n\leq -k,\\
-\sum\limits_{i=1}^{n-k} \left\lfloor \frac{i}{2} \right\rfloor & \textnormal{otherwise},
\end{cases}
\end{align*}
and $f_{n,k}(x,q)$ being a polynomial in $x$ and $q$ satisfying for positive $k$ the recursions
\begin{align}
\label{eq: f quasisym}
f_{n,-k}(x,q)&= (-1)^{n(k+1)}f_{n,k}(x,q),\\
\label{eq: f_0 gleichung 1}
f_{2n,0}(x,q)f_{2n-2,0}(x+2,q)&= -f_{2n-1,1}(x+1,q)^2,\\
\label{eq: f_0 gleichung 2}
f_{2n,0}(x,q)f_{2n,0}(x+2,q)&=f_{2n,1}(x+1,q)^2,
\end{align}
\begin{multline}
\label{eq: f final recursion}
f_{n,k}(x,q)
=\frac{1}{f_{n-2,k}(x+2,q)\left(\frac{n-1}{2}\right)^{[n \in 2\N+1]}} \\
\times \left((x+n)^{[n \in \N\setminus(2\N+k+1)]}\left(q(x+n+1)\right)^{[n \in 2\N+k+2]}f_{n-1,k}(x,q)f_{n-1,k}(x+2,q) \right.\\
-\left. \vphantom{(x+n)^{[n \in \N\setminus(2\N+k+1)]}} (x+1) f_{n-1,k-1}(x+1,q)f_{n-1,k+1}(x+1,q) \right),
\end{multline}
where $[\textnormal{statement}]$ is the \emph{Iverson bracket} which is defined as $1$ if the statement is true and $0$ otherwise, and where we use the convention $0 \in \mathbb{N}$.
%with starting values
%\begin{align*}
%f_{1,k}(x,q)&=\frac{1-(-q)^k}{1+q},\\
%f_{2,k}(x,q)&=\left(\frac{1-(-q)^k}{1+q}\right)^2+(-q)^{k-1}(1+x).
%\end{align*}
\end{thm}

\begin{proof}
Since $d_{n,k}(x,q)$ is a polynomial in $x$ and a Laurent polynomial in $q$ we can write $d_{n,k}(x,q)$ as in \eqref{eq: general form} where $f_{n,k}(x,q)$ is a rational function in $x$ and $q$. Lemma \ref{lem: divisibility by p(n,k)} and Lemma \ref{lem: transponierte matrix} imply that $\frac{d_{n,k}(x,q)}{p_{n,k}(x)}$ is a polynomial in $x$.
Using the $(-q)$-analog of an integer, see \eqref{eq: -q-analog} and the Leibniz formula, we rewrite $d_{n,k}(x,q)$ as
\begin{equation}
\label{eq: d mit leibniz 2}
d_{n,k}(x,q) = \sum_{\sigma \in \mathfrak{S}_n}\prod_{i =1}^n \binom{x+i+\sigma(i)-2}{\sigma(i)-1}[k+\sigma(i)-i]_{-q}.
\end{equation}
Let $\sigma \in \mathfrak{S}_n$.  The exponent of the smallest power of $q$ that appears in the summand associated to $\sigma$ in \eqref{eq: d mit leibniz 2} is 
\begin{equation}
\label{eq: power of q in summand}
\sum_{i: k+\sigma(i)-i<0}(k+\sigma(i)-i).
\end{equation}
It is an easy proof for the reader to show that the minimum of \eqref{eq: power of q in summand} for all $\sigma \in \mathfrak{S}_n$ is exactly $c(n,k)$. Hence $q^{-c(n,k)}d_{n,k}(x,q)$ is a polynomial in $q$, which implies that $f_{n,k}(x,q)$ is a polynomial in $x$ and $q$.\\

It remains to prove that $f_{n,k}$ satisfies the equations \eqref{eq: f quasisym} -- \eqref{eq: f final recursion}, where \eqref{eq: f quasisym} is a direct consequence of Lemma \ref{lem: transponierte matrix}. Equation \eqref{eq: f_0 gleichung 1} and \eqref{eq: f_0 gleichung 2} follow from \eqref{eq: condensed equation} by using Corollary \ref{cor: k=0 and odd n} and \eqref{eq: f quasisym}.
Now let $k \geq 1$. Using the definition of $p_{n,k}(x)$, we can rewrite \eqref{eq: condensed equation} as
\begin{multline*}
d_{n,k}(x,q)q^{c(n-2,k)}\prod_{i=1}^{n-3} \frac{\left\lfloor \frac{i}{2}\right\rfloor!}{i!}\prod_{i=0}^{\left\lfloor\frac{n-k-3}{2}\right\rfloor}(x+k+2i+3)f_{n-2,k}(x+2,q)\\
=\prod_{i=1}^{n-2} \left(\frac{\left\lfloor \frac{i}{2}\right\rfloor!}{i!}\right)^2 \frac{1}{(n-1)}\\
\times \left((x+n)q^{2c(n-1,k)}\prod_{i=0}^{\left\lfloor\frac{n-k-2}{2}\right\rfloor}(x+k+2i+1)(x+k+2i+3)f_{n-1,k}(x,q)f_{n-1,k}(x+2,q) \right.\\
-(x+1)q^{c(n-1,k-1)+c(n-1,k+1)}\prod_{i=0}^{\left\lfloor\frac{n-k-1}{2}\right\rfloor}(x+k+2i+1)\\
\times \left.
\prod_{i=0}^{\left\lfloor\frac{n-k-3}{2}\right\rfloor}(x+k+2i+3)f_{n-1,k-1}(x+1,q)f_{n-1,k+1}(x+1,q)
\right).
\end{multline*}
After cancellation we obtain
\begin{multline*}
d_{n,k}(x,q)f_{n-2,k}(x+2,q)
=\frac{q^{c(n,k)} p_{n,k}(x)}{\left(\frac{n-1}{2}\right)^{[n \in 2\N+1]}}\\
\times \left((x+n)^{[n \in \N\setminus(2\N+k+1)]}(x+n+1)^{[n \in 2\N+k+2]}q^{[n \in 2\N+k+2]}f_{n-1,k}(x,q)f_{n-1,k}(x+2,q) \right.\\
-\left. (x+1) f_{n-1,k-1}(x+1,q)f_{n-1,k+1}(x+1,q) \right).
\end{multline*}
Replacing $d_{n,k}(x,q)$ by $q^{c(n,k)} p_{n,k}(x) f_{n,k}(x,q)$ and further cancellation implies the last recursion \eqref{eq: f final recursion}.
\end{proof}

Computer experiments suggest that $f_{n,k}(x,q)$ is a polynomial with integer coefficients and that the leading coefficient  is either $1$  or $-1$, where we order the monomials $x^a q^b$ with respect to the lexicographic order of $(a,b)$. While we are not able to prove this, we provide the following related statement.

\begin{prop}
The leading coefficient of $f_{n,1}(x,q)$ is $1$,  where we order the monomials $x^a q^b$ with respect to the reverse lexicographic order of $(a,b)$.
\end{prop}

\begin{proof}
First, we calculate the leading coefficient of the highest power of $q$ in $d_{n,1}(x,q)$. This is a polynomial in $x$ of which we calculate the leading coefficient.
We will need the determinantal evaluations
\begin{align}
\label{eq: andere det evaluation}
\nonumber \det_{1 \leq i,j \leq n} \left( \binom{x+i+j-2}{a+j-1} \right) &= \prod_{j=1}^n \frac{(j-1)!}{(a+j-1)!}\prod_{i=1}^{a}(x+j-i),\\
\det_{1 \leq i,j \leq n} \left((-1)^{i+j} \binom{x+i+j-2}{a+j-1} \right) &= \prod_{j=1}^n \frac{(j-1)!}{(a+j-1)!}\prod_{i=1}^{a}(x+j-i),
\end{align}
which can be proved by factoring out common factors in rows and columns and using the Vandermonde determinant evaluation. 
Using the $(-q)$-analog, the determinant $d_{n,1}(x,q)$ can be written as
\begin{equation}
\label{eq: d mit leibniz 3}
d_{n,1}(x,q) = \sum_{\sigma \in \mathfrak{S}_n}\prod_{i =1}^n \binom{x+i+\sigma(i)-2}{\sigma(i)-1}[1+\sigma(i)-i]_{-q}.
\end{equation}
Let $\sigma \in \mathfrak{S}_n$ and $i_1,\ldots, i_l$ be the rows such that $\sigma(i)-i<0$. The exponent of the highest $q$ power that appears in the summand associated to $\sigma$ in \eqref{eq: d mit leibniz 3} is
\begin{align*}
\label{eq: max q power}
-l+\sum_{i: \sigma(i)-i>0} (\sigma(i)-i)= -l+\sum_{j=1}^l(i_j-\sigma(i_j))
\leq -l+l(n-l).
\end{align*}
It is obvious that there exists a $\sigma \in \mathfrak{S}_n$ such that the above inequality is sharp. The maximal exponent is reached for $l=\frac{n-1}{2}$ if $n \equiv 1 \mod 2$ or in the two cases $l=\frac{n-2}{2}$ or $l=\frac{n}{2}$ for $n \equiv 0 \mod 2$.
First, let $n \equiv 1 \mod 2$ and hence $l=\frac{n-1}{2}$. The maximal $q$ power is reached for all $\sigma \in \mathfrak{S}_n$ with $i_j=n-j+1$ and $\sigma(i_j)\leq l$ for all $1 \leq j \leq l$. Hence the coefficient of the  maximal $q$ power in $d_{n,1}(x,q)$ can be written as
\begin{align*}
\det_{1 \leq i,j \leq l} \left(\binom{x+n-l+i+j-2}{j-1}\right)
\det_{1 \leq i,j \leq n-l}\left((-1)^{i+j+l} \binom{x+i+j+l-2}{l+j-1} \right).
\end{align*}
Equation \eqref{eq: andere det evaluation} implies that the leading coefficient with respect to $x$ in the previous expression is equal to
\[
\prod_{i=1}^{n-l}\frac{(i-1)!}{(l+i-1)!}.
\]
By simple manipulations and using $l=\frac{n-1}{2}$ this transforms to
\[
\prod_{i=1}^{n-1}\frac{\left\lfloor\frac{i}{2}\right\rfloor!}{i!},
\]
which proves the claim for odd $n$.\\

Now let $n=2a$. Then the maximal $q$ power is obtained for $l= a$ or $l=a-1$. Analogously to the above case we can write the coefficient of the highest power of $q$ as
\begin{multline*}
\det_{1 \leq i,j \leq a}\left(
\begin{cases}
0  &i=1 \textnormal{ and } j=a,\\
\binom{x+a+i+j-2}{j-1} & \textnormal{otherwise},
\end{cases} \right)
\det_{1 \leq i,j \leq a}\left((-1)^{i+j+a} \binom{x+i+j+a-2}{a+j-1} \right)
\\+\det_{1 \leq i,j \leq a-1}\binom{x+a+i+j-1}{j-1}
\det_{1 \leq i,j \leq a+1} \left( \begin{cases}
0 \quad & i=a+1 \textnormal{ and } j=1,\\
(-1)^{i+j+a-1} \binom{x+i+j+a-3}{a+j-2} &\textnormal{otherwise},
\end{cases} \right).
\end{multline*}
Using Laplace expansion, this becomes
\begin{multline*}
\det_{1 \leq i,j \leq a} \binom{x+a+i+j-2}{j-1}
\det_{1 \leq i,j \leq a}\left((-1)^{i+j+a} \binom{x+i+j+a-2}{a+j-1} \right)\\
+(-1)^a\binom{x+n-1}{a-1}\det_{1 \leq i,j \leq a-1} \binom{x+a+i+j-1}{j-1}
\det_{1 \leq i,j \leq a}\left((-1)^{i+j+a} \binom{x+i+j+a-2}{a+j-1} \right)\\
+\det_{1 \leq i,j \leq a-1}\binom{x+a+i+j-1}{j-1}
\det_{1 \leq i,j \leq a+1} \left(
(-1)^{i+j+a-1} \binom{x+i+j+a-3}{a+j-2}\right)\\
-(-1)^{a}\binom{x+n-1}{a-1}\det_{1 \leq i,j \leq a-1}\binom{x+a+i+j-1}{j-1}
\det_{1 \leq i,j \leq a} \left(
(-1)^{i+j+a} \binom{x+i+j+a-2}{a+j-1}\right).
\end{multline*}
The second and fourth line cancel each other and the degree of $x$ is by \eqref{eq: andere det evaluation} in the first line $a^2$ and for the third line $a^2-1$. Hence the leading coefficient is
\[
\prod_{i=1}^a \frac{(i-1)!}{(a+i-1)!} = \prod_{i=1}^{n-1} \frac{\left\lfloor \frac{i}{2}\right\rfloor!}{i!},
\]
where the equality is an easy transformation.
\end{proof}

One could extend the above proof to the calculation of the leading coefficient of $f_{n,k}(x,q)$ for arbitrary positive $k$ where we order the monomials $x^a q^b$ with respect to the reverse lexicographic order of $(a,b)$. However the leading coefficient of $f_{n,k}(x,q)$ will not be equal to $\pm 1$ for $k \neq 1$.\\

It seems that the polynomial $f_{n,k}(x,q)$ is  for general $k$ irreducible over $\mathbb{Q}[x,q]$. While this appears to be difficult to prove, Theorem \ref{thm: structural thm} implies that $p_{n,k}(x)$ is maximal in $f_{n,k}(x,q)$, i.e., there exists no non-trivial polynomial $p^\prime(x) \in \mathbb{Q}[x]$ dividing $f_{n,k}(x,q)$. For small values of $k$ on the other side we could prove a factorisation of $f_{n,k}(x,q)$ which was already suggested to exist in \cite[pp. 559]{Fischer18}.

\begin{prop}
The function $f_{n,k}(x,q)$ has in the special case for $k=0,1$ the form
\begin{align*}
f_{n,0}(x,q)&=\begin{cases}
0 \qquad & n \textnormal{ odd},\\
(-1)^{\frac{n}{2}} F_{\frac{n}{2}}(x+1,q)^2 & n \textnormal{ even},
\end{cases}\\
f_{n,1}(x,q) &= F_{\left\lfloor \frac{n+1}{2}\right\rfloor}(x,q)F_{\left\lfloor \frac{n}{2}\right\rfloor}(x+2,q),\\
\end{align*}
where $F_{n}(x,q)$ is a polynomial in $x$ and $q$ over $\Q$.
\end{prop}

\begin{proof}
Corollary \ref{cor: k=0 and odd n} implies $f_{2n+1,0}(x,q)=0$. As a consequence of equation \eqref{eq: f_0 gleichung 1} and $f_{2,0}(x,q)=-1$, we can express $(-1)^n f_{2n,0}(x,q)$ as  $F_n(x+1,q)^2$, where $F_n$ is a polynomial in $x$ and $q$. Together with \eqref{eq: f_0 gleichung 1} and \eqref{eq: f_0 gleichung 2} this implies the claim for $f_{n,1}(x,q)$.
\end{proof}

For $k=2$ computer experiments still suggest a decomposition of $f_{n,k}(x,q)$ into two or three factors, depending on the parity of $n$, whereas for $k \geq 3$ the polynomial $f_{n,k}(x,q)$ seems to be irreducible over $\mathbb{Q}[x,q]$. In order to prove the factorisation for $k=2$ one would need to show that the rational function
\[
\frac{q(x+2n+1)(x+2n+2)F_n(x,q)F_n(x+4,q)-nF_{n+1}(x,q)F_{n-1}(x+4,q)}{F_n(x+2,q)},
\]
is a polynomial in $x$ and $q$. Further if one could guess the resulting polynomial, one would obtain a recursion for $F_{n}(x,q)$.\\

The following corollary is a direct consequence of Theorem \ref{thm: structural thm} and the above proposition.

\begin{cor}
Set $Q=q^{-1}+2+q$ and define $\tilde{p}_n(q)$ as the Laurent polynomial
\begin{align*}
\tilde{p}_{2n}(q) &:= 2^{n-1} q^{-\binom{n}{2}}\prod_{i=1}^{n-1}\frac{i!}{(2i)!} F_n(0,q),\\
\tilde{p}_{2n+1}(q) &:= q^{-\binom{n}{2}}\prod_{i=1}^{n}\frac{i!}{(2i-1)!} F_n(2,q),
\end{align*}
%where $(x)_j:=x(x+1)\cdots (x+j-1)$ is the Pochhammer symbol.
Then the $Q$-enumeration of ASMs $A_n(Q)$ is given by
\begin{align*}
A_{2n}(Q)&=2\tilde{p}_{2n}(q)\tilde{p}_{2n+1}(q),\\
A_{2n+1}(Q)&=\tilde{p}_{2n+1}(q)\tilde{p}_{2n+2}(q).\\
\end{align*}
\end{cor}

It is an easy proof for the reader that the Laurent polynomials $\tilde{p}_{n}(q)$ are actually polynomials in $Q$. The above corollary was actually conjectured in \cite[Conjecture 4]{MillsRobbinsRumsey83} and was first proved in \cite{Kuperberg96}.

%%%%%%%%%%%%%%%%%%%%%%%%%%%%%%%%%%%%%%%%%%%%%%%%%%%%%%%%%%%%%%%%%%%%%%%%%%%%%%%%%%%%%%%%%%%%%
%%%%%%%%%%%%%%%%%%%%%%%%%%%%%%%%%%  0,1,2,3,4 enumerations %%%%%%%%%%%%%%%%%%%%%%%%%%%%%%%%%% 
%%%%%%%%%%%%%%%%%%%%%%%%%%%%%%%%%%%%%%%%%%%%%%%%%%%%%%%%%%%%%%%%%%%%%%%%%%%%%%%%%%%%%%%%%%%%%

\section{The $0$-, $1$-, $2$-, $3$- and $4$-enumeration of ASMs}
\label{sec: andere wurzeln}

In this section, we provide factorisations of the determinant $d_{n,k}(x,q)$ where $q$ is a primitive first, second, third, fourth or sixth root of unity. As a consequence of these factorisations we obtain the known formulas for the $1$-, $2$- and $3$-enumeration of ASMs.
The following table shows the connection between the specialisation of $q$ and the weighted enumeration of ASMs.
\bigskip

\begin{center}
\begin{tabular}{c c c }
 $0$-enumeration: & $q=-1$ &(primitive second root of unity),\\
 $1$-enumeration: & $q=-\frac{1}{2}\pm \frac{\sqrt{3}}{2}i$ &(primitive third root of unity),\\
 $2$-enumeration: & $ q= \pm i$ &(primitive fourth root of unity),\\
 $3$-enumeration: & $q=\frac{1}{2}\pm \frac{\sqrt{3}}{2}i$ &(primitive sixth root of unity),\\
 $4$-enumeration: & $q= 1$ &(primitive first root of unity).\\
\end{tabular}
\end{center}
\bigskip

The factorisations in the following theorems can be proved (except for $q=-1$) by induction on $n$ together with \eqref{eq: condensed equation}, Corollary \ref{cor: k=0 and odd n} and cancellation of terms. The case $q=-1$ however is proved solely by using row manipulations.

\begin{thm}
For $q=-1$ holds
\label{thm: second root}
\[
d_{n,1}(x,-1)=\left(2\left\lfloor\frac{n+1}{2}\right\rfloor-1\right)!!\prod_{i=1}^{\left\lfloor \frac{n}{2}\right\rfloor}(x+2i).
\]
\end{thm}
\begin{proof}
The limit of $d_{n,1}(x,q)$ for $q \rightarrow -1$ is
\[
\lim_{q\rightarrow -1} d_{n,1}(x,q) = \left( \binom{x+i+j-2}{j-1}(j-i+1)\right)_{1\leq i,j \leq n}.
\]
The following identity of matrices can be shown by using a variant of the Chu--Vandermonde identity.

\begin{multline*}
\left((-1)^{i+j} \binom{i-1}{j-1}
\right)_{1\leq i,j \leq n}
\times
\left(\binom{x+i+j-2}{j-1}(j-i+1)
\right)_{1\leq i,j \leq n}\\
=
\left(
\delta_{i,j}-\delta_{i,j+1}\left(\frac{j}{x+j}\right)^{\even(j)}
\right)_{1\leq i,j \leq n} \times
 \left(\binom{x+j-1}{j-i}j^{1-\even(i)}\left((x+j)\right)^{even(i)}
\right)_{1\leq i,j \leq n},
\end{multline*}
where $\even(n)$ is equal to $1$ if $n$ is even and $0$ otherwise.
The closed product formula of $d_{n,1}(x,-1)$ is implied by taking the determinant on both sides.
\end{proof}

For the rest of the section it will be convenient to use the Pochhammer symbol, which is defined as  $(x)_j:=x(x+1)\cdots (x+j-1)$.

\begin{thm}
\label{thm: third root}
Let $q$ be a primitive third root of unity. 
Then
\begin{align*}
d_{n,6k+1}(x,q)=& 2^{\left\lfloor\frac{n}{2}\right\rfloor \left\lfloor\frac{n+1}{2}\right\rfloor} \prod_{i=1}^{\left\lfloor\frac{n+1}{2}\right\rfloor}\frac{(i-1)!}{(n-i)!} 
\prod_{i \geq 0}
\left(\frac{x}{2}+3i+1\right)_{\left\lfloor \frac{n-4i}{2}\right\rfloor}
\left(\frac{x}{2}+3i+3\right)_{\left\lfloor \frac{n-4i-3}{2}\right\rfloor}\\
&\times \prod_{i \geq 0}
\left(\frac{x}{2}+n-i+\frac{1}{2}\right)_{\left\lfloor \frac{n-4i-1}{2}\right\rfloor}
\left(\frac{x}{2}+n-i-\frac{1}{2}\right)_{\left\lfloor \frac{n-4i-2}{2}\right\rfloor},\\
d_{n,6k+2}(x,q)=& 2^{\left\lfloor\frac{n}{2}\right\rfloor \left\lfloor\frac{n+1}{2}\right\rfloor}3^{-\left\lfloor\frac{n}{2}\right\rfloor}(1-q)^{n} \prod_{i=1}^{\left\lfloor\frac{n+1}{2}\right\rfloor}\frac{(i-1)!}{(n-i)!}\\
&\times \prod_{i \geq 0}
\left(\frac{x}{2}+n-i\right)_{\left\lfloor\frac{n-4i}{2}\right\rfloor}
\left(\frac{x}{2}+n-i\right)_{\left\lfloor\frac{n-4i-3}{2}\right\rfloor}\\
&\times \prod_{i \geq 0}
\left(\frac{x}{2}+3i+\frac{3}{2}\right)_{\left\lfloor\frac{n-4i-1}{2}\right\rfloor}
\left(\frac{x}{2}+3i+\frac{5}{2}\right)_{\left\lfloor\frac{n-4i-2}{2}\right\rfloor},\\
d_{n,6k+3}(x,q)=&(-q)^n2^{\left\lfloor\frac{n+1}{2}\right\rfloor \left\lfloor\frac{n+2}{2}\right\rfloor-\left\lfloor\frac{n}{2}\right\rfloor} \prod_{i=1}^{\left\lfloor\frac{n+1}{2}\right\rfloor}\frac{(i-1)!}{(n-i)!}\\
&\times \prod_{i \geq 0}
\left(\frac{x}{2}+3i+2\right)_{\left\lfloor\frac{n-4i-1}{2}\right\rfloor}
 \left(\frac{x}{2}+3i+2\right)_{\left\lfloor\frac{n-4i-2}{2}\right\rfloor}\\ 
 &\times \prod_{i \geq 0}
\left(\frac{x}{2}+2\left\lfloor\frac{n}{2}\right\rfloor-i+\frac{1}{2}\right)_{\left\lfloor\frac{n-4i}{2}\right\rfloor}
 \left(\frac{x}{2}+2\left\lfloor\frac{n-1}{2}\right\rfloor-i+\frac{3}{2}\right)_{\left\lfloor\frac{n-4i-3}{2}\right\rfloor},\\
d_{n,6k+4}(x,q)=&q^{-\frac{n}{2}}\, d_{n,6k+2}(x,q),\\
d_{n,6k+5}(x,q)=&q^{-n}\, d_{n,6k+1}(x,q),\\
d_{2n+1,6k}(x,q)=&0,\\
d_{2n,6k}(x,q)=&q^{\frac{n}{2}}
2^{n^2} \prod_{i=1}^{n}\frac{(i-1)!}{(2n-i)!} \\
&\times \prod_{i \geq 0}
\left(\frac{x}{2}+3i+\frac{1}{2}\right)_{n-2i}
\left(\frac{x}{2}+3i+\frac{7}{2}\right)_{n-2(i+1)}
\left(\frac{x}{2}+2n-i\right)_{n-2i-1}^2.
\end{align*}
\end{thm}

\begin{cor}
The number $A_n$ of ASMs of size $n$ is given by
 \[
A_n=\prod_{i=0}^{n-1} \frac{(3i+1)!}{(n+i)!}.
\]
\end{cor}
\begin{proof}
Reorder the terms of the product formula of $d_{n,1}(0,q)$, where $q$ is a primitive third root of unity.
\end{proof}

\begin{thm}
 \label{thm: fourth root}
Let $q$ be a primitive fourth root of unity. Then the following holds
\begin{align*}
d_{2n+1,4k}(x,q)&=0,\\
d_{2n,4k}(x,q)&=(2q)^n \prod_{i=1}^{2n-1}\frac{4^{\left\lfloor \frac{i}{2}\right\rfloor}\left\lfloor \frac{i}{2}\right\rfloor!}{i!}
\prod_{i=1}^{\left\lfloor \frac{n}{2} \right\rfloor}\left(\frac{x}{2}+2i+\frac{1}{2} \right)_{2n-4i+1}
\prod_{i=0}^{\left\lfloor \frac{n-1}{2} \right\rfloor}\left(\frac{x}{2}+2i+\frac{1}{2} \right)_{2n-4i-1},\\
d_{n,4k+1}(x,q)&=2^{\left\lfloor \frac{n}{2}\right\rfloor}\prod_{i=1}^{n-1}\frac{4^{\left\lfloor \frac{i}{2}\right\rfloor}\left\lfloor \frac{i}{2}\right\rfloor!}{i!}
\prod_{i=1}^{\left\lfloor \frac{n}{2} \right\rfloor}\left( \frac{x}{2}+i\right)_{n-2i+1},\\
d_{n,4k+2}(x,q)&=\prod_{i=1}^{n-1}\frac{4^{\left\lfloor \frac{i}{2}\right\rfloor}\left\lfloor \frac{i}{2}\right\rfloor!}{i!}
(2q^{-1})^{\frac{n}{2}}\prod_{i=1}^{\left\lfloor \frac{n}{2} \right\rfloor}\left( \frac{x}{2}+i+\frac{1}{2}\right)_{i}\prod_{i=1}^{\left\lfloor \frac{n-1}{2} \right\rfloor}\left( \frac{x}{2}+i+\frac{1}{2}\right)_{i},\\
d_{n,4k+3}(x,q)&=(-q)^n d_{n,4k+1}(x,q).
\end{align*} 
 \end{thm}

\begin{cor}
The $2$-enumeration $A_n(2)$ of ASMs of size $n$ is given by
 \[
A_n(2)=2^{\binom{n}{2}}.
\]
\end{cor}
\begin{proof}
Reorder the terms of the product formula of $d_{n,1}(0,q)$, where $q$ is a primitive fourth root of unity.
\end{proof}

\begin{thm}
\label{thm: sixth root}
Let $q$ be a primitive sixth root of unity and $k$ an integer. Then the following holds
\begin{align*}
&d_{2n+1,3k}(x,q)= 0,\\
 &d_{2n,3k}(x,q)=q^{2n}c(2n)\prod_{i=0}^{n-1}(x+1+3i)\prod_{i=1}^{n-1}(x+3i)_{2(n-i)},\\
 &d_{n,3k+1}(x,q)=c(n) \prod_{i=0}^{\left\lfloor\frac{n-2}{2}\right\rfloor}(x+2+3i)_{n-1-2i},\\
 &d_{n,3k+2}(x,q)=q^{-n}d_{n,3k+1}(x,q),
\end{align*}
with
\begin{align*}
 c(n)=\begin{cases}
         3^{\frac{(n-2)n}{4}} \prod\limits_{i=0}^{n-1}\frac{\left\lfloor\frac{i}{2}\right\rfloor !}{i!} &
         n \textnormal{ is even},\\
         3^{\frac{(n-1)^2}{4}} \prod\limits_{i=0}^{n-1}\frac{\left\lfloor\frac{i}{2}\right\rfloor !}{i!}
         & \textnormal{otherwise}.
        \end{cases}\\
\end{align*}
\end{thm}

\begin{cor}
The $3$-enumeration $A_n(3)$ of ASMs of size $n$ is given by
 \begin{align*}
 A_{2n+1}(3)&= 3^{n(n+1)}\prod_{i=1}^n\frac{(3i-1)!^2}{(n+i)!^2},\\
 A_{2n+2}(3)&= 3^{n(n+2)} \frac{n!}{(3n+2)!}\prod_{i=1}^{n+1}\frac{(3i-1)!^2}{(n+i)!^2}.
 \end{align*}
\end{cor}
\begin{proof}
Reorder the terms of the product formula of $d_{n,1}(0,q)$, where $q$ is a primitive sixth root of unity.
\end{proof}

\begin{thm}
\label{thm: first root}
For $q=1$
\begin{align*}
d_{n,2k+1}(x,1)&=\prod_{i=2}^{\left\lfloor\frac{n}{2}\right\rfloor}(2i-1)^{-(n+1-2i)} \prod_{i=1}^{\left\lfloor \frac{n}{2}\right\rfloor}(x+2i) p_n(x)p_{n-1}(x),\\
d_{2n,2k}(x,1)&=(-1)^{n}\prod_{i=2}^{n}(2i-1)^{-(2n+1-2i)} \prod_{i=1}^{n}(x+2i-1) p_{2n}(x-1)^2,\\
d_{2n+1,2k}(x,1)&=0,
\end{align*}
where $p_n(x)$ is a polynomial in $x$ satisfying the following recursion
\begin{align*}
p_1(x)& =1,\\
p_3(x) &= 2x+5,\\
p_{2n}(x) &= p_{2n-1}(x+2),\\
p_{2n+1}(x) &= \frac{(x+2n+1)(x+2n+2)p_{2n-1}(x)p_{2n-1}(x+4)-(x+1)(x+2)p_{2n-1}(x+2)^2}{2np_{2n-3}(x+4)}.
\end{align*}
\end{thm}

Computer experiments suggest that the polynomials $p_{n}(x)$ are irreducible over $\mathbb{Q}$.

%%%%%%%%%%%%%%%%%%%%%%%%%%%%%%%%%%%%%%%%%%%%%%%%%%%%%%%%%%%%%%%%%%%%%%%%%%%%%%%%%%%%%%%%%%%%%
%%%%%%%%%%%%%%%%%%%%%%%%%%%%%%%%%%%% Alternative Version  %%%%%%%%%%%%%%%%%%%%%%%%%%%%%%%%%%% 
%%%%%%%%%%%%%%%%%%%%%%%%%%%%%%%%%%%%%%%%%%%%%%%%%%%%%%%%%%%%%%%%%%%%%%%%%%%%%%%%%%%%%%%%%%%%%

\section{Connections to further determinants}
\label{sec: other dets}
The aim of this section is to connect the determinant $d_{n,k}(x,q)$ to further known determinants enumerating ASMs or certain classes of plane partitions. In particular we establish a link to determinantal enumeration formulas for \emph{column strict shifted plane partitions}, the determinant  $\det(M_{\textnormal{ASM}})$ in \cite[Eq. (28)]{BehrendDiFrancescoZinnJustin12} as well as to the generalised Andrews determinant in \cite{CiucuEisenkoelblKrattenthalerZare01}. 
We also present in Section \ref{sec: herleitung2} an alternative proof of Theorem \ref{thm: main thm} using results of \cite{BehrendDiFrancescoZinnJustin12} which is based on the six vertex model approach.
 Further related determinantal formulas can be found for example in \cite[Eq. (29), (65), (66), Theorem 1]{BehrendDiFrancescoZinnJustin12}, \cite[Theorem 3.1]{Lalonde02}, \cite[p. 346]{MillsRobbinsRumsey83} as well as in \cite{Andrews79, CiucuKrattenthaler00, Fischer19b, Kuperberg02, Robbins00}.

\subsection{Column strict shifted plane partitions and an alternative proof of Theorem \ref{thm: main thm}}
\label{sec: herleitung2}
%We will relate the determinant $d_{n,1}(0,q)$ to the determinant $\det\left(M_{\textnormal{ASM}}(n,1,q^{-1}+2+q,1) \right)$, which is defined in \cite[Eq. (28)]{BehrendDiFrancescoZinnJustin12}. Proposition 1 of the same paper states that this determinant is equal to the $(q^{-1}+2+q)$-enumeration of ASMs. Hence, this can be seen as an alternative proof of Theorem \ref{thm: main thm} using the six-vertex model approach. As a by-product, we obtain a connection to the weighted enumeration of column strict shifted plane partitions, a generalisation of DPPs, which we are going to define next.\\

Let $k$ be a non-negative integer. A \emph{column strict shifted plane partition} $\pi$ of class $x$, or short CSSPP of class $x$, is an array of positive integers of the form
\[
\begin{array}{ccc ccc ccc}
\pi_{1,1} &&& \cdots & \cdots & \cdots  &&& \pi_{1,\lambda_1}\\
& \pi_{2,2}  && \cdots & \cdots & \cdots && \pi_{2,\lambda_2} \\
&& \ddots  &&&& \iddots \\
&&& \pi_{l,l} & \cdots & \pi_{l,\lambda_l}
\end{array}
\]
such that 
\begin{itemize}
\item $\lambda_1 \geq \ldots \geq \lambda_l$,
\item the rows are weakly decreasing and the columns are strictly decreasing,
\item the first entry in each row exceeds the number of entries in its row by $x$.
\end{itemize}
For an example of a CSSPP of class $2$, see Figure \ref{fig: CSSPPs and paths} (left).
Column strict shifted plane partitions were first defined in \cite{MillsRobbinsRumsey87}; the above definition was taken from \cite{Fischer19b}. CSSPPs of class $x$ are in bijection with $(2-x)$-DPPs, a generalisation of DPPs defined by Andrews \cite{Andrews79, Andrews80}.\\

\begin{figure}
	\begin{center}
		\begin{tikzpicture}
			\begin{scope}[scale=0.6]
				\node at (1,2) {6};
				\node at (2,2) {6};
				\node at (3,2) {5};
				\node at (4,2) {3};
				\node at (2,1) {4};
				\node at (3,1) {1};
				\end{scope}
  			\begin{scope}[scale=0.6, xshift=6cm, yshift=-2cm]
			 	\foreach \x in {1,...,4}
			 		\foreach \y in {1,...,6}
			 			\draw[fill] (\x,\y) circle [radius=0.07];
			 \draw[line width =0.05cm] (1,6) -- (2,6) -- (2,5) -- (3,5) -- (3,3) -- (4,3) -- (4,1);
			 \draw[line width =0.05cm] (1,4) -- (1,1) -- (2,1);
				\node at (.5,6) {6};
				\node at (1.5,6.35) {6};
				\node at (2.5,5.35) {5};
				\node at (3.5,3.35) {3};
				\node at (.5,4) {4};
				\node at (1.5,1.35) {1};
				
%				\node at (1,3) {7};
%				\node at (2,3) {7};
%				\node at (3,3) {7};
%				\node at (4,3) {6};
%				\node at (5,3) {3};
%				\node at (2,2) {6};
%				\node at (3,2) {5};
%				\node at (4,2) {3};
%				\node at (5,2) {1};
%				\node at (3,1) {4};
%				\node at (4,1) {2};	
%  			\end{scope}
%  			\begin{scope}[scale=0.75, xshift=10cm, yshift=-2cm]
%			 	\foreach \x in {1,...,5}
%			 		\foreach \y in {1,...,7}
%			 			\draw[fill] (\x,\y) circle [radius=0.07];
%			 	
%			 	\draw[line width =0.05cm] (1,7) -- (3,7) -- (3,6) -- (4,6) -- (4,3) -- (5,3) -- (5,1);
%			 	\draw[line width =0.05cm] (1,6) -- (1,5) -- (2,5) -- (2,3) -- (3,3) -- (3,1) -- (4,1);
%			 	\draw[line width =0.05cm] (1,4) -- (1,2) -- (2,2) -- (2,1);
%			 	\node at (.5,7) {7};
%			 	\node at (1.5,7.35) {7};
%			 	\node at (2.5,7.35) {7};
%			 	\node at (3.5,6.35) {6};
%			 	\node at (4.5,3.35) {3};
%			 	\node at (.5,6) {6};
%			 	\node at (1.5,5.35) {5};
%			 	\node at (2.5,3.35) {3};
%			 	\node at (3.5,1.35) {1};
%			 	\node at (.5,4) {4};
%			 	\node at (1.5,2.35) {2};
			\end{scope}
			\begin{scope}[scale=0.55, xshift=19cm, yshift=2.5cm]
		\Tiling
		\TopBoundaryColour{1}{-1}{4}{red!40!white}
	\TopBoundaryColour{2}{-2}{3}{red!40!white}
	\RightBoundaryColour{4}{-1}{-1}{red!40!white}
	\RightBoundaryColour{3}{-2}{0}{red!40!white}
		\draw [line width =0.05cm] (0,4) -- (0,-2) -- ({4*cos(30)},{-6+4*sin(30)}) -- ({4*cos(30)},{4*sin(30)}) -- (0,4);
	
		\draw [line width =0.03cm, color = red, dashed] ({1/2*cos(30)},{3+1/2*sin(30)}) -- ({3/2*cos(30)},{5/2+1/2*sin(30)}) -- ({3/2*cos(30)},{3/2+1/2*sin(30)}) -- ({5/2*cos(30)},{1+1/2*sin(30)}) -- ({5/2*cos(30)},{-1+1/2*sin(30)}) -- ({7/2*cos(30)},{-3/2+1/2*sin(30)}) -- ({7/2*cos(30)},{-7/2+1/2*sin(30)});
		\draw [line width =0.03cm, color = red, dashed]  ({1/2*cos(30)},{1+1/2*sin(30)}) -- ({1/2*cos(30)},{-2+1/2*sin(30)}) -- ({3/2*cos(30)},{-5/2+1/2*sin(30)});
		\end{scope}	
		\end{tikzpicture}
	\end{center}
\caption{A CSSPP $\pi$ of class $2$ with $\rho(\pi)=2$, $\mu(\pi)=2$ (left), its representation as non-intersecting lattice paths (middle) and its corresponding cyclically symmetric lozenge tiling of a cored hexagon (right).}
\label{fig: CSSPPs and paths}
\end{figure}

We will consider two refinements of CSSPPs. Given a CSSPP $\pi$ of class $x$, define $\rho(\pi)$ as the number of rows of $\pi$ and $\mu(\pi)$ as the number of entries $\pi_{i,j}$ such that $\pi_{i,j} \leq x+j-i$. For the CSSPP $\pi$ in Figure \ref{fig: CSSPPs and paths} the statistics are $\rho(\pi)=2$ and $\mu(\pi)=2$.

\begin{lem}
\label{lem: enum of CSSPPs}
The weighted enumeration of CSSPPs of class $x$ with respect to the two statistics $\rho,\mu$ is given by the determinant
\begin{equation}
\label{eq: M_n matrix}
\sum_\pi t^{\rho(\pi)}Q^{\mu(\pi)}=
M_n(x,t,Q):= \det_{0 \leq i,j \leq n-1} \left( 
\delta_{i,j}+t \sum_{l \geq 0} \binom{i}{l} \binom{j+x}{l+x} Q^{j-l}
\right),
\end{equation}
where the sum is over all CSSPPs $\pi$ of class $x$ for which the length of the first row is at most $n$.
\end{lem}
Similar determinantal formulas can be found among others in \cite{Andrews79}, \cite[Eq. (29),(65)]{BehrendDiFrancescoZinnJustin12}, \cite[Eq. (1.5)]{CiucuEisenkoelblKrattenthalerZare01} and \cite[Eq. (7.1)]{Fischer19b}.
\begin{proof}
We use a variation of the well-known Lindstr\"om-Gessel-Viennot Theorem, which can be found for example in the proof of Lemma 3.1 in \cite{CiucuKrattenthaler00}.
Let $s_1,\ldots, s_n$ be starting points, $e_1, \ldots, e_n$ end points and denote by $P(i,j)$ the (weighted) count of lattice paths from $s_i$ to $e_j$.
The (weighted) enumeration of non-intersecting lattice paths from $s_{i_1}, \ldots , s_{i_k}$ to $e_{i_1}, \ldots , e_{i_k}$, where we consider all subsets $\{i_1,\ldots,i_k\}$ of $[n]:= \{1, \ldots, n\}$, is given by
\begin{equation}
\det_{1 \leq i , j \leq n} \left(
\delta_{i,j} + P(i,j)
\right).
\end{equation}
Indeed, by using the Leibniz formula and partially expanding the product, we obtain
\begin{align*}
\det_{1 \leq i , j \leq n} \left( \delta_{i,j} + P(i,j) \right)
= \sum_{\sigma \in \Sgrp_n} \sgn{\sigma} \prod_{i=1}^n \left(\delta_{i,\sigma(i)} + P(i,\sigma(i)) \right)
= \sum_{\sigma \in \Sgrp_n} \sgn{\sigma} \sum_{I_\sigma} \prod_{i \in [n] \setminus I_{\sigma}} P(i,\sigma(i)),
\end{align*}
where $I_\sigma$ is the subset of the fixed points of $\sigma$ for which we have chosen $\delta_{i,\sigma(i)}$ in the product. By exchanging the sums, the above becomes
\begin{align*}
\sum_{I \subseteq [n]} \sum_{\sigma \in \Sgrp_{[n] \setminus I}}\sgn(\sigma) \prod_{i \in [n] \setminus I} P(i,\sigma(i))
= \sum_{I \subseteq [n]} \det_{i,j \in I} \left( P(i,j) \right),
\end{align*}
which is by the Lindstr\"om-Gessel-Viennot Theorem equal to our claim.\\

We interpret a CSSPP $\pi$ of class $x$ as a family of $\rho(\pi)$ many non-intersecting lattice paths. The $i$-th row of $\pi$ corresponds thereby to the path from $(1,\pi_{i,i})$  to $(\lambda_i-i+1,1)$ with horizontal steps at height ($y$-coordinate) $\pi_{i,i+1},\ldots, \pi_{i,\lambda_i}$. For an example see Figure \ref{fig: CSSPPs and paths}. The starting points $s_i$ and end points $e_j$ are then given by $s_i=(1,x+i)$ and $e_j=(j,1)$. The contribution of a row of $\pi$ to the $\mu(\pi)$ statistic corresponds to the number of horizontal steps within the last $j+x$ steps of the corresponding lattice path, with starting point $s_i$ and endpoint $e_j$. Hence the weighted count $P(i,j)$ is equal to
\[
t \sum_{l \geq 0} \binom{i}{j-l} \binom{j+x}{l} Q^{l} = 
t \sum_{l \geq 0} \binom{i}{l} \binom{j+x}{l+x} Q^{j-l},
\]
which proves the claim.
\end{proof}

It was shown in \cite[Proposition 1]{BehrendDiFrancescoZinnJustin12} that 
$(1+q)^{-n}M_n(0,q,q^{-1}+2+q)$
is equal to the $(q^{-1}+2+q)$-enumeration of ASMs, with the parameter $\omega$ appearing in \cite[Eq. (28)]{BehrendDiFrancescoZinnJustin12} equal to $\omega=q/(1+q)$. 
Theorem \ref{thm: main thm} then follows by the next lemma by setting $t=q$, $x=0$ and multiplying by $(1+q)^{-n}$.

\begin{lem}
\label{lem: connecting determinants I}
Let $n$ be a positive integer. Then
\[
M_n(x,t,q^{-1}+2+q) = \det_{0 \leq i,j \leq n-1} \left(
\binom{i+j+x}{j}(1+t(-q)^{j-i})
\right).
\]
\end{lem}
\begin{proof}
The proof consists of a sequence of matrix manipulations using the family of lower triangular matrices defined by
\[
L_{n,x}(a,b):= \left(
\binom{i+x}{j+x}a^{i+x}b^{j+x}
\right)_{0 \leq i,j \leq n-1}.
\]
%Analogously to \cite[Eq. (50)--(53)]{BehrendDiFrancescoZinnJustin12}, these matrices satisfy
It is easy to prove, by using the binomial theorem as well as the identity $\binom{i+x}{l+x}\binom{l+x}{j+x} = \binom{i+x}{j+x}\binom{i-j}{l-j}$, that these matrices satisfy
\begin{align}
\label{eq: L property I}
L_{n,0}(a_1,b_1)L_{n,x}(a_2,b_2)^T&=
\left(
\sum_{l \geq 0} \binom{i}{l}\binom{j+x}{l+x} a_1^i a_2^{j+x}b_1^lb_2^{l+x}
\right)_{0 \leq i,j \leq n-1},\\
\label{eq: L property II}
\det \left( L_{n,x}(a,b) \right) &=
\left(a b \right)^{\binom{n}{2}+n x},\\
\label{eq: L property III}
L_{n,x}(a_1,b_1)L_{n,x}(a_2,b_2)&=L_{n,x}\left(a_1(1+a_2 b_1),\frac{a_2b_1b_2}{1+a_2 b_1}\right),
\end{align}
where $A^T$ denotes the transpose of a matrix $A$.
Set $\gamma= q^{\frac{1}{2}}+q^{-\frac{1}{2}}$. We multiply the $i$-th row of the matrix in the determinant $M_n(x,t,\gamma^2)$ by $\gamma^i$ and the $j$-th column by $\gamma^{-j}$ for all $0 \leq i,j \leq n-1$ and obtain
\[
M_n(x,t,\gamma^2) = \det_{0 \leq i,j \leq n-1}
\left(
\delta_{i,j}+t \sum_{l \geq 0} \binom{i}{l}\binom{j+x}{l+x}\gamma^{i+j-2l}
\right).
\]
By \eqref{eq: L property I}, this can be rewritten as
\[
\det_{0 \leq i,j \leq n-1}
\left(
\Id_n + t L_{n,0}\right(\gamma,\gamma^{-1}\left)L_{n,x}\right(\gamma,\gamma^{-1}\left)^T
\right).
\]
Multiplying the matrix in the above determinant from the left by $L_{n,0}(1,-q^{-\frac{1}{2}})$ and from the right by $L_{n,x}(1,-q^{\frac{1}{2}})^T$ and using \eqref{eq: L property II} and \eqref{eq: L property III}, the previous expression becomes
\[
(-1)^{n x} q^{-\frac{n x}{2}} \det \left(
L_{n,0}(1,-q^{-\frac{1}{2}})L_{n,x}(1,-q^{\frac{1}{2}})^T+t
L_{n,0}(-q^{-1},q^{\frac{1}{2}})L_{n,x}(-q,q^{-\frac{1}{2}})^T
\right)
\]
Finally, using \eqref{eq: L property I} and the Chu-Vandermonde identity proves the claim.
\end{proof}

By setting $x=0$ and $t=q$, Lemma \ref{lem: connecting determinants I} implies $M_n(0,q,q^{-1}+2+q)=(1+q)^{n}d_{n,1}(0,q)$. The determinant $M_n(0,t,Q)$ is the weighted count of CSSPPs of class $0$ with the first row having at most $n$ entries which are in bijection with \emph{cyclically symmetric plane partitions in an $n\times n \times n$ cube}; for a detailed account of this bijection see \cite[p. 76]{MillsRobbinsRumsey82}. Hence we obtain a connection between weighted cyclically symmetric plane partitions and weighted ASMs. In particular, we have
\begin{align}
\label{eq: 1-enum}  \sum_{A} 1 &= \zeta_6^{-n} \sum_{\pi} \zeta_3^{\rho(\pi)} ,\\
\label{eq: 2-enum} \sum_{A} 2^{\mu(A)} &= (\sqrt{2} \zeta_8)^{-n} \sum_{\pi} \zeta_4^{\rho(\pi)} 2^{\mu(\pi)},\\
\label{eq: 3-enum}\sum_{A} 3^{\mu(A)} &= (\sqrt{3}\zeta_{12})^{-n}\sum_{\pi} \zeta_6^{\rho(\pi)} 3^{\mu(\pi)},\\
\label{eq: 4-enum} \sum_{A} 4^{\mu(A)} &= 2^{-n}\sum_{\pi}  4^{\mu(\pi)} ,
\end{align}
where the left sum is over all ASMs $A$ of size $n$ and $\mu(A)$ denotes the number of $-1$'s of $A$,  the right sum is over all CSSPPs $\pi$ of class $0$ and entries at most $n$ and $\zeta_l$ denotes the $l$-th primitive root of unity $\zeta_l=e^{\frac{2\pi i}{l}}$.
It would be of interest to have ``bijective'' proofs for these identities, i.e., a usual bijection for  \eqref{eq: 4-enum} and probabilistic bijections for \eqref{eq: 1-enum}--\eqref{eq: 3-enum}.

For $x=0, t=1$ and $q=\zeta_3$, Lemma \ref{lem: connecting determinants I} implies
that the number of cyclically symmetric plane partitions inside a box with side length $n$ is given by
$\zeta_6^n d_{n,3}(0,\zeta_3)$. Together with Theorem \ref{thm: third root} this gives another proof of the weak Macdonald conjecture, which was first proved by Andrews \cite{Andrews79}.

%%%%%%%%%%%%%%%%%%%%%%%%%%%%%%%%%%%%%%%%%%%%%%%%%%%%%%%%%%%%%%%%%%%%%%%%%%%%%%%%%%%%%%%%%%%%%
%%%%%%%%%%%%%%%%%%%%%%%%%%%%%%%%  Connections to Andrews Det %%%%%%%%%%%%%%%%%%%%%%%%%%%%%%%% 
%%%%%%%%%%%%%%%%%%%%%%%%%%%%%%%%%%%%%%%%%%%%%%%%%%%%%%%%%%%%%%%%%%%%%%%%%%%%%%%%%%%%%%%%%%%%%

\subsection{Connections to the Andrews determinant}
\label{sec: verbindungen}

In \cite{Fischer18} Fischer conjectured that $d_{n,1}(x,q)$ is connected to the determinant
\begin{align}
\label{eq: dpp det}
\det_{1 \leq i,j \leq n} \left(\binom{x+i+j-2}{j-1}+q \delta_{i,j}\right).
\end{align}
It was shown by Andrews \cite[Theorem $3^\prime$]{Andrews79} that for $q=1$ the above determinant has a closed product formula and counts the number of $(2-x)$-DPPs with parts less than $n$, i.e., CSSPPs of class $x$ where the first row has at most $n$ entries.
In \cite{CiucuEisenkoelblKrattenthalerZare01} and \cite[Lemma 3.1]{CiucuKrattenthaler00}, it was shown by Ciucu et al.\ that the above determinant is equal to the weighted enumeration of CSSPPs of class $x$, where the first row is at most of length $n$.
Further they proved that the evaluation of the determinant can be expressed by a closed product formula if $q$ is a sixth root of unity.
Comparing the factorisations of $d_{n,k}(x,q)$ and of the above determinant implies
\begin{equation}
\label{eq: d(n,k) and Andrews}
d_{n,3-k}(x,q^2)=q^{-n}  \det_{1 \leq i,j \leq n} \left(\binom{x+i+j-2}{j-1}+q^k \delta_{i,j}\right),
\end{equation}
where $q$ is a primitive sixth root of unity. Alternatively, this can be deduced as follows.\\

First, we note that CSSPPs of class $x$ for which the first row has at most $n$ entries are in bijection with cyclically symmetric lozenge tilings of a cored hexagon with side lengths $n+x,n,n+x,n,n+x,n$ where an equilateral triangle with side length $x$ is removed from the centre. This was first described by Krattenthaler in \cite{Krattenthaler06}. The non-intersecting lattice path description for a CSSPP corresponds thereby to the configuration of lozenges in a fundamental region of the cyclically symmetric lozenge tiling, which is bound by thick black lines in Figure \ref{fig: CSSPPs and paths} (right). For better readability, we do not distinguish between a CSSPP and its corresponding lozenge tiling. Using the lozenge tiling interpretation, we define an involution $\tau$ on CSSPPs of class $x$ and with at most $n$ entries in the first row, by reflecting a lozenge tiling along the horizontal axis. This axis is depictured in Figure \ref{fig: CSSPP involution} by a dotted horizontal line. The involution $\tau$ was first considered for the case $x=2$ in \cite{MillsRobbinsRumsey83} and for arbitrary $x$ in \cite[p. 54]{MillsRobbinsRumsey87} and \cite[p. 1145]{Krattenthaler06}. It has the properties
\begin{equation}
\label{eq: props of tau}
\mu( \pi) = \mu(\tau(\pi)), \qquad \rho(\pi) = n - \rho(\tau(\pi)).
\end{equation}
\begin{figure}
\begin{center}
\begin{tikzpicture}
	\begin{scope}[scale=0.55, xshift=-5cm]
	\draw [line width=0.03cm, dashed] (-8,-1) -- (6,-1);
	\Tiling	
	\TopBoundaryColour{1}{-1}{4}{red!40!white}
	\TopBoundaryColour{2}{-2}{3}{red!40!white}
	\RightBoundaryColour{4}{-1}{-1}{red!40!white}
	\RightBoundaryColour{3}{-2}{0}{red!40!white}
	\draw [line width =0.03cm, color = red, dashed] ({1/2*cos(30)},{3+1/2*sin(30)}) -- ({3/2*cos(30)},{5/2+1/2*sin(30)}) -- ({3/2*cos(30)},{3/2+1/2*sin(30)}) -- ({5/2*cos(30)},{1+1/2*sin(30)}) -- ({5/2*cos(30)},{-1+1/2*sin(30)}) -- ({7/2*cos(30)},{-3/2+1/2*sin(30)}) -- ({7/2*cos(30)},{-7/2+1/2*sin(30)});
	\draw [line width =0.03cm, color = red, dashed]  ({1/2*cos(30)},{1+1/2*sin(30)}) -- ({1/2*cos(30)},{-2+1/2*sin(30)}) -- ({3/2*cos(30)},{-5/2+1/2*sin(30)});
	\TopBoundaryColour{6}{-6}{0}{blue!40!white}
	\TopBoundaryColour{5}{-5}{1}{blue!40!white}
	\LeftBoundaryColour{4}{0}{4}{blue!40!white}
	\LeftBoundaryColour{3}{-1}{3}{blue!40!white}
	\draw [line width =0.03cm, color = blue, dashed] ({1/2*cos(30)},{-11/2+1/2*sin(30)}) -- ({3/2*cos(30)},{-11/2+3/2*sin(30)}) -- ({3/2*cos(30)},{-9/2+3/2*sin(30)}) -- ({5/2*cos(30)},{-9/2+5/2*sin(30)}) -- ({5/2*cos(30)},{-5/2+5/2*sin(30)}) -- ({7/2*cos(30)},{-5/2+7/2*sin(30)}) -- ({7/2*cos(30)},{-1/2+7/2*sin(30)});
	\draw [line width =0.03cm, color = blue, dashed] ({1/2*cos(30)},{-7/2+1/2*sin(30)}) -- ({1/2*cos(30)},{-5/2+1/2*sin(30)}) -- ({3/2*cos(30)},{-5/2+3/2*sin(30)}) -- ({3/2*cos(30)},{-1/2+3/2*sin(30)});
		\draw [line width =0.05cm] (0,0) -- (0,-2) -- ({4*cos(30)},{-6+4*sin(30)}) -- ({4*cos(30)},{4*sin(30)}) -- (0,0);
	\end{scope}
\end{tikzpicture}
\end{center}
\caption{\label{fig: CSSPP involution}The cyclically symmetric lozenge tiling of a cored hexagon of Figure \ref{fig: CSSPPs and paths} together with its representation as non-intersecting lattice paths (red) and the construction of the involution $\tau$ using non-intersecting lattice paths (blue).}
\end{figure}
Indeed, the non-intersecting lattice path of the lozenge tiling (drawn in red in Figure \ref{fig: CSSPP involution}) is mapped under $\tau$ to the lattice path in blue of the same figure. By the definition of $\mu$, the statistic $\mu(\pi)$ is exactly the number of horizontal lozenges within the trapezoidal region (bounded by the thick black lines in Figure \ref{fig: CSSPP involution}) and therefore invariant under reflection along the horizontal axis. This implies the first equation in \eqref{eq: props of tau}. The second equation follows by the fact, that there are exactly $n$ lozenges within the trapezoidal region touching the top boundary. Each of these lozenges is either passed by a red path--- there are $\rho(\pi)$ many--- or the starting point of a blue path--- there are $\rho(\tau(\pi))$ many.

As a consequence of \eqref{eq: props of tau} we obtain the identity
\begin{equation}
\label{eq: t reversing identity}
M_n(x,t,Q)=t^n M_n(x,t^{-1},Q).
\end{equation}
Equation \eqref{eq:  d(n,k) and Andrews} follows from Lemma \ref{lem: connecting determinants I} and \eqref{eq: t reversing identity} by setting $t=\zeta_6^k$ and $q=\zeta_3$.\\

The identity \eqref{eq:  d(n,k) and Andrews} implies that all factorisations of the determinant in \eqref{eq: dpp det} are covered by the determinant $d_{n,k}(x,q)$, where $q$ is a third root of unity. This is in particular interesting since the determinant $d_{n,k}(x,q)$ can be evaluated by using the Desnanot-Jacobi identity (Condensation method) which is ``easier'' than the previously known method for proving the factorisations of the determinant in \eqref{eq: dpp det}, see \cite{CiucuEisenkoelblKrattenthalerZare01}.\\

Very recently it was shown by Fischer \cite{Fischer19b} that for $q=1$ the determinant in \eqref{eq: dpp det} also enumerates $(n,x)$-alternating sign trapezoids, which generalise alternating sign triangles (ASTs). This complements the equinumerousity of DPPs and ASTs (and hence ASMs) by the equinumerousity of a ``one parameter generalisation'' of DPPs and ASTs. An interpretation of $x$ in $d_{n,k}(x,q)$, e.g. in the form of a one parameter generalisation of ASMs, would therefore be very interesting and could give important insight into the nature of the equinumerousity between ASMs, DPPs and ASTs. %Such possible interpretations of $x$  will be studied in a forthcoming work.

\section*{Acknowledgements}

The author is indebted to the careful review and the helpful suggestions of two anonymous referees. In particular, the connection to the determinant $\det(M_{\ASM})$ in \cite[Eq. (28)]{BehrendDiFrancescoZinnJustin12} by using Lemma \ref{lem: connecting determinants I} and a direct proof  for \eqref{eq: d(n,k) and Andrews} by using Lemma \ref{lem: connecting determinants I} and the involution $\tau$ was communicated to the author by one of the referees.
The author also wants to thank Christoph Koutschan for helpful comments.

\begin{appendix}
\section{An overview of enumeration formulas connected to $d_{n,k}(x,q)$}
\label{app: Data}
It turns out that some specialisations of $x,q,k$ in $d_{n,k}(x,q)$ are known enumeration formulas.
In the following we list those specialisations we found. In particular there is always a combinatorial interpretation for $d_{n,k}(0,q)$ where $k$ is an integer and $q$ is a sixth root of unity but not equal to $-1$. All formulas can be proved by using induction on $n$. We denote by  $\zeta_l$ the $l$-th root of unity $\zeta_l= e^{\frac{2 \pi i}{l}}$ and further use the notation $A_{QT}^{(1)}(4n,x)$ of \cite{Kuperberg02}.

\begin{equation}
\zeta_6^{-n} d_{2n,0}(0,\zeta_3) = \left( A_{QT}^{(1)}(4n,1)\right)^2= \left(\#(\textnormal{ASMs of size }n) \right)^4,
\end{equation}
\begin{equation}
\zeta_4^{-n} d_{2n,0}(0,\zeta_4) = \left( A_{QT}^{(1)}(4n,2)\right)^2,
\end{equation}
\begin{equation}
\zeta_3^{-n} d_{2n,0}(0,\zeta_6) = \left( A_{QT}^{(1)}(4n,3)\right)^2 = 3^{n(n-1)}\left(\#(\textnormal{ASMs of size }n) \right)^2,
\end{equation}

\bigskip

\begin{equation}
d_{n,1}(0,-1) = (\textnormal{0-enumeration of ASMs of order }n ) = n!,
\end{equation}
\begin{equation}
d_{n,1}(0,\zeta_3) = \zeta_6^{n-1} d_{n-1,3}(2,\zeta_3) = \#(\textnormal{ASMs of size }n),
\end{equation}
\begin{align}
\nonumber d_{n,1}(0,\zeta_4) &= (\textnormal{2-enumeration of ASMs of order }n )\\
\nonumber &= \#(\textnormal{perfect matchings of an order }n \textnormal{ Atzec diamond} )\\
&= \#(\textnormal{Gelfand-Tsetlin patterns with bottom row }1,2,\ldots,n) =2^{\binom{n}{2}},
\end{align}
\begin{equation}
d_{n,1}(0,\zeta_6) = (\textnormal{3-enumeration of ASMs of size }n),
\end{equation}

\bigskip

\begin{equation}
\zeta_{12}^n d_{n,2}(0,\zeta_3) = \sqrt{3}^{[n \equiv 1 \mod 2]}  \left(\#(\textnormal{half turn symmetric ASMs of size }n)\right)^2,
\end{equation}
\begin{align}
\nonumber \zeta_{4}^n (1-\zeta_4)^n d_{n,2}(0,\zeta_4) =&  \textnormal{Total dimension of the homology of a free }\\
&2\textnormal{-step nilpotent Lie algebra of rank } n\footnotemark,
\end{align}\footnotetext{See \cite[Theorem 1.1]{GrassbergerKingTirao02}. This sequence further coincides with tilings of a half-hexagon of side length $n$ with glued sides, see \cite[Eq. (3.5)]{DiFrancescoZinnJustinZuber05}.}

\bigskip

\begin{align}
\nonumber \zeta_6^{n} d_{n,3}(0,\zeta_3) &= \#(\textnormal{cyclic symmetric plane partitions in an } n-\textnormal{cube}) \\
 &= \#(\textnormal{half turn symmetric ASMs of size } 2n)/\#(\textnormal{ASMs of size }n),
\end{align}

\bigskip

\begin{align}
\zeta_6^{n-1} d_{2(n-1),2}(1,\zeta_3) &=\frac{\zeta_6^{n+1}}{\sqrt{-3}} d_{2n-1,2}(-1,\zeta_3) = \#(\textnormal{ASMs with two U-turn sides of size } 4n),
\end{align}
\begin{equation}
\frac{\zeta_6^{n+1}}{\sqrt{-3}} d_{2n-1,2}(-2,\zeta_3) = \#(\textnormal{quarter turn symmetric ASMs of size }4n),
\end{equation}
\begin{equation}
\zeta_4^{n}d_{2n,2}(1,\zeta_4)= \frac{\zeta_8^{2n+1}}{\sqrt{2}}d_{2n+1,2}(-1,\zeta_4)= 4^{n^2},
\end{equation}
\begin{equation}
 \frac{\zeta_8^{2n-1}}{\sqrt{2}}d_{2n-1,2}(1,\zeta_4)= 4^{n(n+1)}.
\end{equation}
\end{appendix}

\bibliographystyle{abbrv}
\bibliography{LiteraturListe}

\end{document}